\documentclass[11pt]{article}

\usepackage{amsfonts,amssymb,amstext}
\usepackage{amsmath}
\usepackage{amstext}
\usepackage{amsthm}

\makeatletter

\newcommand{\beq}{\begin{equation}}
\newcommand{\eeq}{\end{equation}}
\newcommand{\bea}{\begin{eqnarray}}
\newcommand{\eea}{\end{eqnarray}}
\newcommand{\beaa}{\begin{eqnarray*}}
\newcommand{\eeaa}{\end{eqnarray*}}

\newcommand{\n}{\noindent}
\newcommand{\q}{\quad}
\newcommand{\qq}{\qquad}

\newtheorem{thm}{Theorem}[section]

\newtheorem{pp}{Proposition}[section]

\newtheorem{exa}{Example}[section]
\newtheorem{rem}{Remark}[section]

\newtheorem{asn}{Assumption}[section]

\newcommand{\g}{\gamma}

\newcommand{\I}{\varphi}
\newcommand{\G}{\Gamma}
\newcommand{\de}{\delta}
\newcommand{\De}{\Delta}

\newcommand{\al}{\alpha}
\newcommand{\la}{\lambda}
\newcommand{\f}{\infty}
\newcommand{\vs}{\varepsilon}
\newcommand{\cd}{\cdot}
\newcommand{\si}{\sigma}

\newcommand{\be}{\beta}
\newcommand{\ta}{\theta}

\newcommand{\sm}{\setminus}

\newcommand {\ol} {\overline}

\newcommand {\s} {\section}
\newcommand {\sn} {\subsection}
\newcommand {\ssn} {\subsubsection}

\newcommand{\p}{\partial}

\newcommand{\tr}{\hbox{tr}}

\newcommand{\ConvD}{\overset{d}{\rightarrow}}
\newcommand{\ConvP}{\overset{P}{\rightarrow}}

\newcommand{\Var}{\mathrm{Var}}
\newcommand{\E}{\mathbb{E}}

\newcommand{\wh}{\widehat}
\newcommand {\uu}{{\bf u}}
\newcommand {\vv}{{\bf v}}

\newcommand{\best}{\widehat{{\ta}}_T}

\numberwithin{equation}{section}

\usepackage{color}

\begin{document}

\title{Statistical inference for %L\'evy-driven
stationary linear models with tapered data}

\author{Mamikon S. Ginovyan\\
\small Department of Mathematics and Statistics\\
\small Boston University, Boston, MA, USA\\
%111 Cummington Mall, Boston, MA 02215, USA
\small e-mail: ginovyan@math.bu.edu\\
\\
Artur A. Sahakyan\\
\small Department of Mathematics and Mechanics\\
\small Yerevan State University, Yerevan, Armenia\\
\small e-mail: sart@ysu.am}

\date{}
%\date{\today}
\maketitle

\begin{abstract}
\noindent
In this paper, we survey some recent results on statistical inference
(parametric and nonparametric statistical estimation, hypotheses testing)
about the spectrum of stationary models with tapered data, as well as,
a question concerning robustness of inferences, carried out on a linear
stationary process contaminated by a small trend.%, to this departure from stationarity.
We also discuss some question concerning tapered Toeplitz matrices and operators,
central limit theorems for tapered Toeplitz type quadratic functionals,
and tapered Fej\'er-type kernels and singular integrals.
These are the main tools for obtaining
the corresponding results, and also are of interest in themselves.
The processes considered will be discrete-time and continuous-time Gaussian,
linear or L\'evy-driven linear processes with memory.

\end{abstract}

\maketitle

\tableofcontents

\section{Introduction}
\label{INT}

\vskip2mm
\noindent
Let $\{X(t), \ t\in \mathbb{U}\}$ be a centered real-valued stationary process
with spectral density $f(\la)$, $\la\in \Lambda$, and covariance function
$r(t)$, $t\in \mathbb{U}$. We consider simultaneously the
continuous-time (c.t.)\ case, where $\mathbb{U}=\mathbb{R}:=(-\f,\f)$, and the
discrete-time (d.t.)\ case, where $\mathbb{U}=\mathbb{Z}:=\{0,\pm 1,\pm 2,\ldots \}$.
The domain $\Lambda$ of the frequency variable $\la$ is $\Lambda=\mathbb{R}$
in the c.t.\ case, and $\Lambda:= [-\pi.\pi]$ in the d.t.\ case.

We want to make statistical inferences (parametric and nonparametric estimation,
and hypotheses testing) about the spectrum of $X(t)$. In the classical setting,
the inferences are based on an observed finite realization $\mathbf{X}_T$
of the process $X(t)$:
$\mathbf{X}_T:=\{X(t), \, t\in D_T\},$
where $D_T:=[0, T]$ in the c.t.\ case and $D_T:=\{1,\ldots,T\}$
in the d.t.\ case.
%$$
%\mathbf{X}_T:=\{X(t), \,\, 0\le t\le T  \, \, \text{(in the c.t.\ case), or} \,\,
%t={1,\ldots,T} \, \, \text{(in the d.t.\ case)}\}.
%$$

A sufficiently developed inferential theory is now available for stationary models
based on the standard (non-tapered) data $\mathbf{X}_T$. %For instance,
We cite merely the following references Anh et al. \cite{ALS2}, Avram et al. \cite{AvLS},
Casas and Gao \cite{CG}, Dahlhaus \cite{D}, Dahlhaus and Wefelmeyer \cite{DW},
Dzhaparidze \cite{Dz1,Dz},  Dzhaparidze and Yaglom \cite{DY}, Fox and Taqqu \cite{FT1}, Gao \cite{Go},
Gao et al. \cite{GAHT}, Ginovyan  \cite{G1988a, G1988b, G2003a, G2011a,G2003g,G2018g,G2020g},
Giraitis et al. \cite{GKSu}, Giraitis and Surgailis \cite{GSu}, Guyon \cite{Gu}, Hannan \cite{H},
Has'minskii and Ibragimov \cite{IH-2}, Heyde and Dai \cite{HD}, Heyde and Gay \cite{HG},
Ibragimov \cite{I1963,I1967}, Ibragimov and Khas'minskii \cite{IH-3},
Leonenko and Sakhno \cite{LS}, Millar \cite{M-1},
Osidze \cite{O1,O2}, Taniguchi \cite{Tan87}, Taniguchi and Kakizawa \cite{TK}, Tsai and Chan \cite{TC},
Walker \cite{Wal}, Whittle \cite{W},
where can also be found additional references.

In the statistical analysis of stationary processes, however, the data are frequently
tapered before calculating the statistic of interest, and the statistical inference
procedure, instead of the original data $\mathbf{X}_T$, is based on the {\it tapered data}:
$\mathbf{X}^h_T:=\{h_T(t)X(t), \, t\in D_T\}\},$
where $D_T:=[0, T]$ in the c.t.\ case and $D_T:=\{1,\ldots,T\}$
in the d.t.\ case, and $h_T(t):= h(t/T)$ with $h(t)$, $t\in\mathbb{R}$ being a {\it taper function}.
%$$
%\mathbf{X}^h_T:=\{h_T(t)X(t), \,\, 0\le t\le T  \, \, \text{(in the c.t.\ case), or} \,\,
%t={1,\ldots,T} \, \, \text{(in the d.t.\ case)}\},
%$$
%where $h_T(t):= h(t/T)$ with $h(t)$, $t\in\mathbb{R}$ being a {\it taper function}.

The use of data tapers in nonparametric time series was suggested by Tukey \cite{Tu}.
The benefits of tapering the data have been widely reported in the literature
(see, e.g., Brillinger \cite{Bri2},  Dahlhaus \cite{D1}--\cite{D4}, \cite{D5},
Dahlhaus and K\"unsch \cite{DK}, Guyon \cite{Gu}, and references therein).
For example, data-tapers are introduced to reduce the so-called 'leakage effects',
that is, to obtain better estimation of the spectrum of the model in the case
where it contains high peaks.
Other application of data-tapers is in situations in which some of the
data values are missing. Also, the use of tapers leads to bias reduction,
which is especially important when dealing with spatial data. In this case,
the tapers can be used to fight the so-called 'edge effects'.

In this paper, we survey some recent results on statistical inference
(parametric and nonparametric statistical estimation, and hypotheses testing)
about the spectrum of stationary models with tapered data, as well as,
a question concerning robustness of inferences, carried out on a linear
stationary process contaminated by a small trend. %, to this departure from stationarity.
We also discuss some questions concerning tapered Toeplitz matrices and operators,
central limit theorems for tapered Toeplitz type quadratic functionals,
and tapered Fej\'er-type kernels and singular integrals.
These are the main tools for obtaining
the corresponding results, and also are of interest in themselves.
The processes considered will be discrete-time and continuous-time Gaussian,
linear or L\'evy-driven linear processes with memory.

\paragraph{Some notation and conventions.}
The following notation and conventions are used throughout the paper.\\
The symbol '$:=$' stands for 'by definition';
c.t.: = continuous-time; d.t.:= discrete-time;
s.d.:= spectral density; c.f.:= covariance function;
CLT:= central limit theorem.
The symbols '$\ConvP$' and '$\ConvD$' stand for convergence in probability and in distribution, respectively.
The notation $X_T\ConvD \eta \sim N(0,\sigma^2)$ as $T\to\f$
will mean that the distribution of the random variable
$X_T$ tends (as $T\to\f$) to the centered normal
distribution with variance $\sigma^2$.
$\E[\cdot]$: = expectation operator;
$\tr[A]$: = trace of an operator (matrix) $A$;
${\mathbb I}_A(\cdot)$: = indicator of a set $A\subset \Lambda$;
WN$(0,1)$: = standard white-noise.
The standard symbols $\mathbb{N}$, $\mathbb{Z}$ and $\mathbb{R}$
denote the sets of natural, integer and real numbers, respectively;
$\mathbb{N}_0:=\mathbb{N}\cup0$.
By $\Lambda$ we denote the frequency domain, that is, $\Lambda:=\mathbb{R}$
in the c.t.\ case, and $\Lambda:= [-\pi.\pi]$ in the d.t.\ case.
By $L^p(\mu):=L^p(\Lambda,\mu)$ ($p\geq $1) we denote the weighted
Lebesgue space with respect to the measure $\mu$, and by $||\cdot||_{p,\mu}$
we denote the norm in $L^p(\mu)$.
In the special case where $d\mu(\la)=d\la$, we will use the notation
$L^p$ and $||\cdot||_{p}$, respectively.
The letters $C$ and $c$ with or without indices
are used to denote positive constants, the values of which can vary from
line to line.
Also, in the d.t.\  case all the considered functions are assumed to be
$2\pi$-periodic and periodically extended to $\mathbb{R}$.
%Also, we will assume that $f, g \in L^1(\mathbb{R})$, and with no
%loss of generality, that $g\ge 0$.

\paragraph{The structure of the paper.}
The rest of the paper is structured as follows.
In Section \ref{Pre} we specify the model of interest - a stationary process,
recall some key notions and results from the theory of stationary processes,
and introduce the data tapers and tapered periodogram. % based on the tapered data.
In Section \ref{NP} we discuss the nonparametric estimation problem.
We analyze the asymptotic properties, involving asymptotic unbiasedness,
bias rate convergence, consistency, a central limit theorem and asymptotic
normality of the empirical spectral functionals.
%, needed to prove consistency and asymptotic normality of the Whittle estimator.
In Section \ref{PW} we discuss the parametric estimation problem.
We present sufficient conditions for consistency and asymptotic normality
of minimum contrast estimator based on the Whittle contrast functional
for stationary linear models with tapered data.
Section \ref{GF} is devoted to the construction of goodness-of-fit
tests for testing hypotheses that the hypothetical spectral density
of a stationary Gaussian model has the specified form, based on
the tapered data.
A question concerning robustness of inferences, carried out on a linear
stationary process contaminated by a small trend is discussed in Section \ref{Rob}.
In Section \ref{methods} we briefly discuss the methods and tools, used to prove
the results stated in Sections \ref{NP}--\ref{Rob}.

\section{Preliminaries}
\label{Pre}
In this section we specify the model of interest - a stationary process,
and introduce the data tapers and tapered periodogram.

\sn{The model}
%\paragraph{Second-order (wide-sense) stationary process.}
\ssn{Second-order (wide-sense) stationary process}
Let $\{X(u), \ u\in \mathbb{U}\}$ be a centered real-valued second-order (wide-sense)
stationary process defined on a probability space $(\Omega, \mathcal{F}, P)$ with covariance
function $r(t)$, that is, $\E[X(u)]=0$, $r(u)=\E[X(t+u)X(t)]$, $u, t\in\mathbb{U}$,
%\[ \E[X(u)]=0, \q r(u)=\E[X(t+u)X(t)], \q u, t\in\mathbb{U}, \]
where $\E[\cdot]$ stands for the expectation operator with respect to measure $P$.
We consider simultaneously the c.t.\ case, where
$\mathbb{U}=\mathbb{R}:=(-\f,\f)$, and the d.t.\ case, where
$\mathbb{U}=\mathbb{Z}:=\{0,\pm 1,\pm 2,\ldots \}$.
We assume that $X(u)$ is a {\it non-degenerate process}, that is,
${\rm Var}[X(u)]=\E|X(u)|^2=r(0)>0$.
(Without loss of generality, we assume that $r(0)=1$).
In the c.t.\ case the process $X(u)$ is also assumed
mean-square continuous, that is, $\mathbb{E}[X(t)-X(s)]^2\to0$ as
$t\to s$.
This assumption is equivalent to that of the covariance function $r(u)$
be continuous at $u=0$ (see, e.g., Cram\'er and Leadbetter \cite{CrL}, Section 5.2).

By the Herglotz theorem in the d.t.\ case, and the Bochner-Khintchine
theorem in the c.t.\ case (see, e.g., Cram\'er and Leadbetter \cite{CrL}),
there is a finite measure $\mu$ on $(\Lambda, \mathfrak{B}(\Lambda))$,
where $\Lambda=\mathbb{R}$ in the c.t.\ case, and $\Lambda= [-\pi.\pi]$
in the d.t.\ case, and $\mathfrak{B}(\Lambda)$ is the Borel $\si$-algebra on $\Lambda$,
such that for any $u\in \mathbb{U}$ the covariance function $r(u)$
admits the following {\sl spectral representation}:
\beq
\label{i1}
r(u)=\int_\Lambda \exp\{i\la u\}d\mu(\la), \q u\in\mathbb{U}.
\eeq
\n

The measure $\mu$ in (\ref{i1}) is called the {\sl spectral measure} of the
process $X(u)$. The function  $F(\la): =\mu[-\pi,\la]$ in the d.t.\ case
and $F(\la): =\mu[-\f,\la]$ in the c.t.\ case, is called the {\sl spectral
function} of the process $X(t)$.
If $F(\la)$ is absolutely continuous (with respect to Lebesgue measure),
then the function $f(\la):=dF(\la)/d\la$ is called the {\sl spectral density}
of the process $X(t)$. Notice that if the spectral density $f(\la)$ exists, then
$f(\la)\geq 0$, $f(\la)\in L^1(\Lambda)$, and \eqref{i1} becomes
\begin{equation}
\label{mo1}
r(u)=\int_\Lambda\exp\{i\la u\}f(\la)d\la, \q  u\in\mathbb{U}.
\end{equation}
Thus, the covariance function $r(u)$ and the spectral
function $F(\la)$ (resp. the spectral density function $f(\la)$) are equivalent
specifications of the second order properties for a stationary process
$\{X(u), \ u\in\mathbb{U}\}$.

%\paragraph{Linear processes. Existence of spectral density functions.}
\ssn{Linear processes. Existence of spectral density functions}
We will consider here stationary processes possessing spectral density
functions. For the following results we refer to
Cram\'er and Leadbetter \cite{CrL}, Doob \cite{Doob}, and Ibragimov and Linnik \cite{IL}.
\begin{itemize}
\item[(a)]
The spectral function $F(\la)$ of a d.t.\ stationary process $\{X(u), \,u \in \mathbb{Z}\}$
is absolutely continuous (with respect to the Lebesgue measure), $F(\la)=\int_{-\pi}^\la f(x)dx$,
if and only if it can be represented as an infinite moving average:
\begin{equation}
\label{dlp}
X(u) = \sum_{k=-\f}^{\f}a(u-k)\xi(k), \qq
         \sum_{k=-\f}^{\f}|a(k)|^2 < \f,
\end{equation}
where $\{\xi(k), k\in \mathbb{Z}\}\sim$ WN(0,1) is a standard white-noise, that is,
a sequence of orthonormal random variables.
\item[(b)]
The covariance function $r(u)$ and the spectral density $f(\la)$ of $X(u)$ are given by formulas:
\begin{equation}
\label{dcv}
r(u)= \E X(u)X(0)=\sum_{k=-\f}^{\f}a(u+k) a(k),
\end{equation}
and
\beq
\label{dsd}
f(\la) = \frac{1}{2\pi} |\wh a(\la)|^2 =
\frac{1}{2\pi}\left|\sum_{k=-\f}^{\f}a(k)e^{-ik\la}\right|^2,
\quad \la\in[-\pi,\pi].
\eeq
\item[(c)]
In the case where $\xi(k)$ is a sequence of Gaussian random variables,
the process $X(u)$ is Gaussian.
\end{itemize}
Similar results hold for c.t.\ processes. Indeed, the following holds.
\begin{itemize}
\item[(a)]
The spectral function $F(\la)$ of a c.t.\ stationary process $\{X(u), \,u \in \mathbb{R}\}$
is absolutely continuous (with respect to Lebesgue measure), $F(\la)=\int_{-\f}^\la f(x)dx$,
if and only if it can be represented as an infinite continuous moving average:
\begin{equation}
\label{clp}
X(u)=\int_{\mathbb{R}} a(u-t)d\xi(t),, \qq
         \int_{\mathbb{R}}|a(t)|^2dt < \f,
\end{equation}
where $\{\xi(t), t\in \mathbb{R}\}$ is a process with orthogonal increments and
$\E|d\,\xi(t)|^2 = dt$.
\item[(b)]
The covariance function $r(u)$ and the spectral density
$f(\la)$ of $X(u)$ are given by formulas:
\begin{equation}
\label{ccv}
r(u)= \E X(u)X(0)=\int_{\mathbb{R}} a(u+x)a(x)dx,
\end{equation}
and
\beq
\label{csd}
f(\la) = \frac{1}{2\pi} |\wh a(\la)|^2 =
\frac{1}{2\pi}\left|\int_{\mathbb{R}} e^{-i\la t}a(t)dt\right|^2,
\quad \la\in\mathbb{R}.
\eeq
\item[(c)]
In the case where $\xi(t)$ is a Gaussian process,
the process $X(u)$ is Gaussian.
\end{itemize}

%\paragraph{L\'evy-driven linear process.}
\ssn{L\'evy-driven linear process}
We first recall that a L\'evy process, $\{\xi(t),\ t\in\mathbb{R}\}$ is a process with
independent and stationary increments, continuous in probability, with sample-paths which
are right-continuous with left limits (c$\grave{a}$dl$\grave{a}$g) and
$\xi(0) = \xi(0-) = 0$.
The Wiener process $\{B(t), \ t\ge 0\}$ and the centered Poisson process
$\{N(t)-\E N(t), \ t\ge 0\}$ are typical examples of centered L\'evy processes.
A L\'evy-driven linear process $\{X(t),\ t\in\mathbb{R}\}$ is a real-valued c.t.\
stationary process defined by \eqref {clp},
%\begin{equation}\label{lp}
%X(t)=\int_{\mathbb{R}} a(t-s)\xi(ds),
%\end{equation}
%where  $a(\cdot)$ is a function from $L^2(\mathbb{R})$, and $\xi(t)$
where $\xi(t)$ is a L\'evy process satisfying the conditions:
$\E \xi(t)=0$, $\E \xi^2(1)=1$ and $\E\xi^4(1)<\infty$.
In the case where $\xi(t)=B(t)$, $X(t)$ is a Gaussian process
 (see, e.g., Bai et al. \cite{BGT2}):
%\[
%\text{$\E \xi(t)=0$, $\E \xi^2(1)=1$ and $\E\xi^4(1)<\infty$.}
%\]

%In the case where $\xi(t)=B(t)$, $X(t)$ is a Gaussian process.

The function $a(\cdot)$ in representations \eqref{dlp} and \eqref{clp}
plays the role of a {\it time-invariant filter}, and the linear processes
defined by \eqref{dlp} and \eqref{clp} can be viewed as the output of
a linear filter $a(\cdot)$ applied to the process $\xi(t)$,
called the innovation or driving process of $X(t)$.
%According to Wold's fundamental theorem, a very broad class of stationary
%processes can be represented as linear processes.

Processes of the form \eqref{dlp} and \eqref{clp} appear in many fields of science
(economics, finance, physics, etc.), and cover large classes of popular models
in time series modeling.
For instance, the classical autoregressive moving average (ARMA) models
and their continuous counterparts the c.t.\ autoregressive moving average (CARMA) models
are of the form \eqref{dlp} and \eqref{clp}, respectively,
and play a central role in the representations of stationary time series
(see, e.g.,  Brockwell \cite{Br1}, Brockwell and Davis \cite{BD}).

%\paragraph{Dependence (memory) structure of the model}
\ssn{Dependence (memory) structure of the model.}
In the frequency domain setting, the statistical and spectral analysis
of stationary processes requires
{\em two types of conditions\/} on the spectral density $f(\la).$
The first type controls the {\em singularities} of $f(\la)$,
and involves the {\em dependence (or memory) structure } of the
process, while the second type -- controls the {\em smoothness} of $f(\la).$
The memory structure of a stationary process is essentially a
measure of the dependence between all the variables in the process,
considering the effect of all correlations simultaneously.
Traditionally memory structure has been defined in the time domain
in terms of decay rates of the autocorrelations, or in the
frequency domain in terms of rates of explosion of low frequency spectra
(see, e.g., Beran et al. \cite{BFGK}, Giraitis et al. \cite{GKSu}, Gu\'egan \cite{Gu1}).
It is convenient to characterize the memory structure in terms
of the spectral density function.

We will distinguish the following types of stationary models:

(a) short memory (or short-range dependent),

(b) long memory (or long-range dependent),

(c) intermediate memory (or anti-persistent).

\vskip2mm
\noindent
{\em Short-memory models.}
Much of statistical inference
is concerned with {\it short-memory} stationary models,
where the spectral density $f(\lambda)$ of the
model is bounded away from zero and infinity,
that is, there are constants $C_1$ and $C_2$ such that
\begin{equation*}
0< C_1 \le f(\la) \le C_2 <\f.
\end{equation*}

A typical d.t.\ short memory model example is the stationary
Autoregressive Moving Average (ARMA)$(p,q)$ process $X(t)$
defined to be a stationary solution of the difference equation:
\begin{equation*}
\psi_p(B)X(t)=\theta_q(B)\varepsilon(t), \q t\in\mathbb{Z},
\end{equation*}
where $\psi_p$ and $\theta_q$ are polynomials of degrees $p$ and $q$,
respectively, $B$ is the backshift operator defined by $BX(t)=X(t-1)$,
and $\{\vs(t), t\in\mathbb{Z}\}$ is a d.t.\ white noise,
that is, a sequence of zero-mean, uncorrelated random variables
with variance $\si^2$.
%that is, $E[\vs(t)]=0$ and $E[\vs(t)\vs(s)]=\si^2\de(t,s)$,
The covariance $r(k)$ of (ARMA)$(p,q)$ process is
exponentially bounded:
\[
|r(k)|\le Cr^{-k}, \q k=1,2,\ldots;\q0<C<\f; \,\, 0<r<1,
\]
and the spectral density $f(\la)$ is a rational function
(see, e.g., Brockwell and Davis \cite{BD}, Section 3.1):
%
%e9 ###
\begin{equation}
\label{arma}
f(\la) = \frac{\si^2}{2\pi}\cd\frac
{|\theta_q(e^{-i\la})|^2}{|\psi_p(e^{-i\la})|^2}.
\end{equation}

A typical c.t.\ short-memory model example is the stationary
c.t.\ ARMA$(p,q)$ processes, denoted by CARMA$(p,q)$,
The spectral density function $f(\la)$ of a CARMA$(p,q)$ process $X(t)$ is given by
the following formula (see, e.g., Brockwell \cite{Br1}):
\begin{equation}
\label{CAS}
f(\la)=\frac{\si^2}{2\pi}\cd\frac{|\be_q(i\la)|^2}{|\al_p(i\la)|^2},
\end{equation}
where $\al_p(z)=z^p-\al_pz^{(p-1)}-\cdots-\al_1$ and $\be_q(z)=1+\be_1z+\cdots+\be_qz^q$ are
polynomials of degrees $p$ and $q$, respectively.

Another important c.t.\ short-memory model is the
{\it Ornstein-Uhlenbeck} process, which is a Gaussian stationary
process with covariance function $r(t)=\si^2e^{-\al |t|}$ ($t\in\mathbb{R}$),
and spectral density
\begin{equation}
\label{OU}
f(\la) = \frac{\si^2}{\pi}\cd\frac {\al^2}{\la^2+\al^2},\q \al>0,  \, \la\in\mathbb{R}.
\end{equation}

\vskip2mm
\noindent
{\em Discrete-time long-memory and anti-persistent models.}
Data in many fields of science (economics, finance, hydrology, etc.),
however, is well modeled by stationary processes whose spectral densities are
{\it unbounded} or {\it vanishing} at some fixed points  (see, e.g.,
Beran et al. \cite{BFGK}, Gu\'egan \cite{Gu1}, and references therein).

A {\it long-memory} model is defined to be a
stationary process with {\it unbounded} spectral density,
and an {\it anti-persistent} model -- a stationary
process with {\it vanishing} (at some fixed points) spectral density.

In the discrete context, a basic long-memory model
is the Autoregressive Fractionally Integrated Moving Average
(ARFIMA)$(0,d,0))$ process $X(t)$ defined to be a stationary
solution of the difference equation
(see, e.g., Brockwell and Davis \cite{BD}, Section 13.2):
\begin{equation*}
(1-B)^dX(t)=\varepsilon(t), \q0<d<1/2,
\end{equation*}
where $B$ is the backshift operator and
$\varepsilon(t)$ is a d.t.\ white noise defined above.
The spectral density $f(\la)$ of $X(t)$ is given by
%
%e10 ###
\begin{equation}
\label{MT0}
f(\la)=|1-e^{-i\la}|^{-2d}=(2\sin(\la/2))^{-2d},
\q0<\la\le\pi, \q0<d<1/2.
\end{equation}
Notice that $f(\la)\thicksim c\, |\la|^{-2d}$ as $\la\to0$, that is, $f(\la)$ blows up
at $\la=0$ like a power function, which is the typical behavior of
a long memory model.

A typical example of an {\it anti-persistent} model is the
ARFIMA$(0,d,0)$ process $X(t)$ with spectral density specified by (\ref{MT0})
%$f(\la)=|1-e^{-i\la}|^{-2d}$
with $d<0$, which vanishes at $\la=0$.
Note that the condition $d<1/2$ ensures that
$\int_{-\pi}^\pi f(\la)d\la<\f$,
implying that the process $X(t)$ with spectral density (\ref{MT0}) is well defined because
$E|X(t)|^2=\int_{-\pi}^\pi f(\la)d\la.$

Data can also occur in the form of a realization of a 'mixed'
short-long-intermediate-memory stationary process $X(t)$.
%with spectral density
%$$f(\la)=f_I(\la)f_L(\la)f_S(\la),$$
%where $f_I(\la)$, $f_L(\la)$ and $f_S(\la)$ are the intermediate,
%long- and short-memory components, respectively.
A well-known example of such a process, which appears in many applied
problems, is an ARFIMA$(p,d,q)$ process $X(t)$ defined to be a
stationary solution of the difference equation:
\begin{equation*}
\psi_p(B)(1-B)^dX(t)=\theta_q(B)\varepsilon(t), \q d<1/2,
\end{equation*}
where $B$ is the backshift operator, $\varepsilon(t)$ is a
d.t.\ white noise, and $\psi_p$ and $\theta_q$ are
polynomials of
degrees $p$ and $q$, respectively.
The spectral density $f_X(\la)$ of $X(t)$ is given by
%
%e11 ###
\begin{equation}
\label{AA}
f_X(\la)=|1-e^{-i\la}|^{-2d}f(\la), \q d<1/2,
\end{equation}
where $f(\la)$ is the spectral density of an ARMA$(p,q)$ process,
given by (\ref{arma}).
\noindent
Observe that for $ 0<d<1/2$ the model $X(t)$ specified by the spectral density (\ref{AA})
displays long-memory, for $d<0$ -- intermediate-memory,
and for $d=0$ -- short-memory.
For $d \ge1/2$ the function $f_X(\la)$ in (\ref{AA}) is not integrable,
and thus it cannot represent a spectral density of a stationary
process. Also, if $d \le-1$, then the series $X(t)$ is not
invertible in the sense that it cannot be used to recover a white noise
$\varepsilon(t)$ by passing $X(t)$ through a linear filter
(see, e.g., Brockweel and Davis \cite{BD}).

Another important long-memory model is the fractional Gaussian noise (fGn).
To define the fGn we first introduce the {\it fractional Brownian motion\/} (fBm)
\mbox{$\{B_H(t), t\in\mathbb{R}\}$} with Hurst index $H$, $0<H<1$,
defined to be a centered Gaussian $H$-self-similar process having stationary
increments (see, e.g., Samorodnisky and Taqqu \cite{ST}).
Then the increment process $\{X(k):= B_H(k+1)-B_H(k), k\in\mathbb{Z}\},$
% \[\{X(k):= B_H(k+1)-B_H(k), k\in\mathbb{Z}\},\]
called {\it fractional Gaussian noise\/} (fGn), is a d.t.\
centered Gaussian stationary process with spectral density function:
%
%e13 ###
\begin{equation}
\label{MT2}
f(\la)=c\, |1-e^{-i\la}|^{2}\sum_{k=-\f}^\f|\la+2\pi k|^{-(2H+1)},
\q-\pi\le\la\le\pi,
\end{equation}
where $c$ is a positive constant.

It follows from (\ref{MT2}) that $f(\la)\thicksim c\, |\la|^{1-2H}$
as $\la\to0$, that is, $f(\la)$ blows up if $H>1/2$
and tends to zero if $H<1/2$.
Also, comparing (\ref{MT0}) and (\ref{MT2}), we observe that,
up to a constant, the spectral density of fGn has the same
behavior at the origin as ARFIMA$(0,d,0)$ with $d=H-1/2.$

Thus, the fGn $\{X(k), k\in\mathbb{Z}\}$ has long-memory
if $1/2<H<1$ and is anti-percipient if $0<H<1/2$.
The variables $X(k)$, $k\in\mathbb{Z}$, are independent if $H=1/2$.
For more details we refer to Samorodnisky and Taqqu \cite{ST}.

\vskip2mm
\noindent
{\em Continuous-time long-memory and anti-persistent models.}
In the continuous context, a basic process which has commonly been used
to model long-range dependence is the fractional Brownian motion (fBm)
$B_H$ with Hurst index $H$, defined above, which can be regarded as a
Gaussian process having a 'spectral density':
%of the form
%
%e14 ###
\begin{equation}
\label{lr3}
f(\la)=c|\la|^{-(2H+1)}, \q c>0, \ \,\,
0<H<1, \  \,\, \la\in\mathbb{R}.
\end{equation}
The form (\ref{lr3}) can be understood in a generalized sense
(see, e.g., Yaglom \cite{Ya}),
since the fBm $B_H$ is a nonstationary process.
%(see, also, Anh et al. \cite{ALM} and Gao et al. \cite{GAHT}).

A proper stationary model in lieu of fBm is the
{\it fractional Riesz-Bessel motion} (fRBm),
introduced in Anh et al. \cite{AAR}, and defined as a
c.t.\ Gaussian process $X(t)$ with
spectral density
\begin{equation}
\label{rb1}
f(\la)=c\,|\la|^{-2\al}(1+\la^2)^{-\be},\q\la\in\mathbb{R},
\, 0<c<\f, \, 0<\al<1, \,\be>0.
\end{equation}
The exponent $\al$ determines the long-range dependence,
while the exponent $\be$
indicates the second-order intermittency of the process
(see, e.g., Anh et al. \cite{ALM} and Gao et al. \cite{GAHT}).

Notice that the process $X(t)$, specified by the spectral density (\ref{rb1}),
is stationary if $0<\al<1/2$
and is non-stationary with stationary increments if $1/2\le\al<1.$
\noindent Observe also that the function $f(\la)$ in (\ref{rb1}) behaves
as $O(|\la|^{-2\al})$ as $|\la|\to0$ and as $O(|\la|^{-2(\al+\be)})$
as $|\la|\to\f$.
Thus, under the conditions $0<\al<1/2$, $\be>0$ and $\al+\be>1/2,$
the function $f(\la)$ in (\ref{rb1}) is well-defined for both
$|\la|\to0$ and $|\la|\to\f$
due to the presence of the component $(1+\la^2)^{-\be}$, $\be>0$,
which is the Fourier transform of the Bessel potential.

Comparing (\ref{lr3}) and (\ref{rb1}), we observe that the
spectral density of fBm is the limiting case as $\be\to0$
that of fRBm with Hurst index $H=\al-1/2.$

Another important c.t.\ long-memory model is the CARFIMA$(p,H,q)$ process.
The spectral density function $f(\la)$ of a CARFIMA$(p,H,q)$ process $X(t)$ is given by
the following formula (see, e.g., Brockwell \cite{Br1}, and Tsai and Chan \cite{TC}):
\begin{equation}
\label{CAA}
f(\la)=\frac{\si^2}{2\pi}\G(2H+1)\sin(\pi H)|\la|^{1-2H}\frac{|\be_q(i\la)|^2}{|\al_p(i\la)|^2},
\end{equation}
where $\al_p(z)=z^p-\al_pz^{(p-1)}-\cdots-\al_1$ and $\be_q(z)=1+\be_1z+\cdots+\be_qz^q$ are
polynomials of degrees $p$ and $q$, respectively.
Notice that for $H=1/2$, the spectral density given by \eqref{CAA} becomes
that of the short-memory CARMA$(p,q)$ process, given by \eqref{CAS}.

\sn{Data tapers and tapered periodogram}

Our inference procedures will be based on the tapered data $\mathbf{X}_T^h$:
\beq \label{t2}
\mathbf{X}_T^h: =
\left \{
\begin{array}{ll}
 \{h_T(t)X(t),  \, t=1,\ldots, T\} & \mbox{in the d.t.\ case},\\
\{h_T(t)X(t),  \, 0\leq t\leq T\} & \mbox{in the c.t.\ case},
\end{array}
\right.
\eeq
where
\beq
\label{Tap}
h_T(t):= h(t/T)
\eeq
with $h(t)$, $t\in\mathbb{R}$ being a {\it taper function}. % to be specified below.

Throughout the paper, we will assume that the taper function  $h(\cdot)$
satisfies the following assumption.

\begin{asn}
\label{(T)}
%{\bf Assumption (T).}
{\rm The taper $h:\mathbb{R}\to\mathbb{R}$ is a continuous nonnegative
function of bounded variation and of bounded support $[0, 1]$, such that
$H_k\neq 0$, where}
\beq
\label{t4}
H_{k}: = \int_0^1h^k(t)dt, \q k\in \mathbb{N}:=\{1,2,\ldots\}.
\eeq
\end{asn}

\n
{\it Note.} The case $h(t)={\mathbb I}_{[0,1]}(t)$, where ${\mathbb I}_{[0,1]}(\cdot)$
denotes the indicator of the segment $[0,1]$, will be referred to as
the {\it non-tapered} case.
\begin{rem}
\label{TH}
{\rm For the d.t.\ case, an example of a taper function $h(t)$ satisfying Assumption \ref{(T)}
is the Tukey-Hanning taper function $h(t)=0.5(1-\cos(\pi t))$ for $t\in [0, 1]$.
For the c.t.\ case, a simple example of a taper function $h(t)$ satisfying
Assumption \ref{(T)} is the function $h(t)=1-t$ for $t\in [0, 1]$.
More examples of taper functions satisfying Assumption \ref{(T)}
can be found in Dahlhaus \cite{D5} and in Guyon \cite{Gu}.}
\end{rem}

Denote by $H_{k,T}(\la)$ the {\it tapered Dirichlet type kernel}, defined by
\beq \label{t3}
H_{k,T}(\la): =
\left \{
\begin{array}{ll}
\sum_{t=1}^T h_{T}^k(t)e^{-i\la t} & \mbox{in the d.t.\ case},\\
\\[-1mm]
\int_0^T h_{T}^k(t)e^{-i\la t}dt & \mbox{in the c.t.\ case}.
\end{array}
\right.
\eeq

Define the finite Fourier transform of the tapered data (\ref{t2}):
\beq
\label{t5}
d^h_T(\la): =
\left \{
\begin{array}{ll}
\sum_{t=0}^T h_{T}(t)X(t)e^{-i\la t} & \mbox{in the d.t.\ case},\\
\\[-1mm]
 \int_0^Th_T(t)X(t)e^{-i\la t}dt & \mbox{in the c.t.\ case}.
\end{array}
\right.
\eeq
and the tapered periodogram $I^h_{T}(\lambda)$ of the process $X(t)$:
\bea
\label{t6}
I^h_{T}(\lambda):= \frac1{C_T}\,d^h_T(\la)d^h_T(-\la)=
\left \{
\begin{array}{ll}
\frac 1{C_T}\left |\sum_{t=0}^T h_{T}(t)X(t)e^{-i\la t}\right|^2 &
\mbox{in the d.t.\ case},\\
\\[-1mm]
\frac1{C_T}\left|\int_0^Th_T(t)X(t)e^{-i\la t}dt\right|^2 &
\mbox{in the c.t.\ case}.
\end{array}
\right.
\eea
where
\beq
\label{t55}
C_T:= 2\pi H_{2,T}(0)=2\pi\int_0^Th_T^2(t)dt=2\pi H_2\,T \neq 0.
\eeq
Notice that for non-tapered case ($h(t)={\mathbb I}_{[0,1]}(t)$),
we have $C_T= 2\pi T$.

%\s {Asymptotic properties of empirical spectral functionals}%{Estimation of spectral functionals.}
%\label{NP}

\section{Nonparametric Estimation Problem}
\label{NP}

Suppose we observe a finite realization $\mathbf{X}_T:=\{X(u)$, $0\le u\le T$
(or $u={1,\ldots,T}$ in the d.t.\ case)\}
of a centered stationary  process $X(u)$ with an {\it unknown\/} spectral
density function $f(\lambda)$, $\lambda \in \Lambda$.
We assume that $f(\lambda)$ belongs to a given (infinite-dimensional)
class $\mathcal{F} \subset L^p:=L^p(\Lambda)$ $(p \ge 1)$
of spectral densities possessing some specified smoothness properties.
The problem is to estimate the value  $J(f)$ of a given
functional $J(\cdot)$ at an {\it unknown} 'point' $f \in \mathcal{F}$
on the basis of an observation $\mathbf{X}_T,$ and investigate the asymptotic
(as $T \to \infty$) properties of the suggested estimators, depending on the
dependence structure of the model $X(u)$ and the smoothness structure of the
'parametric' set $\mathcal{F} \subset L^p(\Lambda)$ $(p \ge 1)$.

Linear and non-linear functionals of the periodogram play a key
role in the parametric estimation of the spectrum of stationary
processes, when using the minimum contrast estimation method with
various contrast functionals (see, e.g., Anh et al. \cite{ALS2},
Dzhaparidze \cite{Dz}, Guyon \cite{Gu}, Leonenko and Sakhno \cite{LS},
Taniguchi and Kakizawa \cite{TK}, and references therein).
In this section, we review the asymptotic properties, involving
asymptotic unbiasedness, bias rate convergence, consistency,
a central limit theorem and asymptotic normality of the empirical
spectral functionals based on the tapered data. Some of these properties were
discussed and proved in Ginovyan and Sahakyan \cite{GS2019, GS2019a}.
For non-tapered case, these properties were established in the papers
Ginovyan \cite{G1994, G2011a}.
The results stated in this section are used to prove consistency and
asymptotic normality of the minimum contrast estimator
based on the Whittle contrast functional for stationary linear models
with tapered data (see Section \ref{PW}).
Here we follow the papers Ginovyan \cite{G1995,G2011a,G2011b},
and Ginovyan and Sahakyan \cite{GS2019, GS2019a}.

\sn {Estimation of linear spectral functionals}
\label{Lem}
We are interested in the nonparametric estimation problem, based on the
tapered data \eqref{t2}, of the following linear spectral functional:
\beq
\label{t1}
J=J(f,g): = \int_\Lambda f(\lambda)g(\lambda)d\lambda,
\eeq
where $g(\lambda) \in L^q(\Lambda)$, \,$1/p + 1/q = 1$.

As an estimator $J_T^h$ for functional $J(f)$, given by \eqref{t1},
based on the tapered data \eqref{t2}, we consider the averaged tapered
periodogram (or a simple 'plug-in' statistic), defined by
\bea
\label{t7}
J_T^h &=& J(I^h_{T}):= \int_\Lambda I^h_{T}(\la)g(\lambda)d\lambda,
\eea
where $I^h_{T}(\lambda)$ is the tapered periodogram of the process
$X(t)$ given by \eqref{t6}.
Denote
\beq
\label{tq-1}
Q_T^h:=
 \left \{
 \begin{array}{ll}
\sum_{t=1}^T\sum_{s=1}^T\widehat g(t-s)h_T(t)h_T(s)X(t)X(s)
& \mbox{in the d.t.\ case},\\
\\
\int_0^T\int_0^T \widehat g(t-s)h_T(t)h_T(s)X(t)X(s)\,dt\,ds
& \mbox{in the c.t.\ case},
\end{array}
\right. \eeq
where $\widehat g(t)$ is the Fourier transform of function $g(\la)$:
\beq
\label{Ft7}
\widehat g(t):=\int_{\Lambda} e^{i\lambda t} g(\lambda) d\lambda,
\quad t\in\Lambda.
\eeq
In view of (\ref{t6}) and (\ref{t7}) -- (\ref{Ft7}) we have
\bea
\label{t77}
J_T^h =C_T^{-1}Q_T^h,
%=\left \{
%\begin{array}{ll}
%\frac 1{C_T}\sum_{t=1}^T \sum_{s=1}^T\widehat g(t-s)h_T(t)h_T(s)X(t)X(s) &
%\mbox{in the d.t.\ case},\\
%\\[-1mm]
%\frac 1{C_T}\int_0^T\int_0^T \widehat g(t-s)h_T(t)h_T(s)X(t)X(s)\,dt\,ds &
%\mbox{in the c.t.\ case},
%\end{array}
%\right.
\eea
where $C_T$ is as in (\ref{t55}).
We will refer to $g(\la)$ and to its Fourier transform $\widehat g(t)$ as a
{\it generating function} and {\it generating kernel} for the functional
$J_T^h$, respectively.

Thus, to study the asymptotic properties of the estimator $J_T^h$, we have
to study the asymptotic distribution (as $T\to\f$) of the tapered Toeplitz type
quadratic functional $Q_T^h$ given by (\ref{tq-1}) (for details see Section \ref{CLT}).

\sn{Asymptotic unbiasedness}
We begin with the following assumption.

\begin{asn}
\label{(A1)}
%{\bf Assumption (A1).}
{\rm The function
\bea
\label{au1}
\Psi(u) = \int_{\Lambda} f(v)g (u + v)\,dv
\eea
belongs to $L^1(\Lambda)\cap L^2(\Lambda)$ and is continuous at $u = 0$.}
\end{asn}

\begin{thm}
\label{AU1}
Let the functionals $J:=J(f,g)$ and $J_T^h: = J(I^h_{T},g)$ be defined
by (\ref {t1}) and (\ref{t7}), respectively.
Then under Assumptions \ref{(T)} and \ref{(A1)}
the statistic $J_T^h$ is an asymptotically unbiased estimator for $J(f)$,
that is, the following relation holds:
\bea
\label{au2}
\lim_{T \to \f}[E (J_T^h) - J] = 0.
\eea
\end{thm}

\begin{rem}\label{r1}
{\rm Using H\"older inequality, it can easily be shown that if
$f\in L^1(\Lambda)\cap L^{p_1}(\Lambda)$
and $g\in L^1(\Lambda)\cap L^{p_2}(\Lambda)$ with $1\le p_1,p_2\le\f$,
$1/p_1+1/p_2\le 1$, then the relation \eqref{au2} is satisfied.}
\end{rem}

Under additional smoothness conditions on functions $f(\la)$ and $g(\la)$ we can
estimate the rate of convergence in \eqref{au2}.
To state the corresponding result, we first introduce
some notation and assumptions.

Given numbers $p\ge 1$, $0<\al<1$, $r\in \mathbb{N}_0:=\mathbb{N}\cup\{0\}$,
where $\mathbb{N}$ is the set of natural numbers, we set $\be=\al+r$
and denote by $H_p(\be)$ the $L^p$-H\"older class, that is,
the class of those functions $\psi(\la)\in L^p(\Lambda)$, which have  $r$-th
derivatives in $L^p(\Lambda)$ and with some positive constant $C$ satisfy
$$||\psi^{(r)}(\cdot+h)-\psi^{(r)}(\cdot)||_p\le C|h|^\al.$$
\begin{asn}
\label{(A2')}
%{\bf Assumption (A2').}
{\rm We say that a pair of integrable functions $(f(\la), g(\la))$, $\la\in\Lambda$,
satisfies condition $(\mathcal{H})$, and write $(f, g)\in (\mathcal{H})$,
if $f\in H_p(\be_1)$ for $\be_1>0$, $p>1$ and $g\in H_q(\be_2)$ for $\be_2>0$,
$q>1$ with $1/p+ 1/q = 1$, and one of the conditions a) -- d)
is satisfied:

\n a) $\be_1>1/p$,\,\, $\be_2>1/q$,

\n b) $\be_1 \le 1/p$,\,\, $\be_2\le 1/q$ and $\be_1+\beta_2>1/2$,

\n c) $\be_1>1/p$,\,\, $1/q-1/2<\be_2 \le 1/q$,

\n d) $\be_2>1/q$,\,\, $1/p-1/2 <\be_1 \le 1/p$.}
\end{asn}
\begin{rem}
\label{RH1}
{\rm In Ginovian \cite{G1994} it was proved that if $(f, g)\in (\mathcal{H})$,
then there exist numbers
$p_1$ $(p_1>p)$ and $q_1$ $(q_1>q)$, such that $H_p(\be_1)\subset L_{p_1}$,
$H_q(\be_2)\subset L_{q_1}$ and $1/{p_1}+1/{q_1}\le1/2$.}
\end{rem}

\begin{asn}
\label{(A2)}
%{\bf Assumption (A2).}
{\rm The spectral density $f$ and the generating function $g$ are such that
$f, g\in L^1(\Lambda)\cap L^2(\Lambda)$ and $g$ is of bounded variation.}
\end{asn}

The following theorem controls the bias $E(J_T^h)-  J$ and provides sufficient conditions
assuring the proper rate of convergence of bias to zero, necessary for asymptotic
normality of the estimator $J_T^h$. Specifically, we have the following result.
%the proof of which is given in Ginovyan and Sahakyan \cite{GS2019}.

\begin{thm}
\label{AU3}
Let the functionals $J:=J(f,g)$ and $J_T^h: = J(I^h_{T},g)$ be defined
by (\ref {t1}) and (\ref{t7}), respectively.
Then under Assumptions \ref{(T)} and \ref{(A2')} (or \ref{(A2)}),
the following asymptotic relation holds:
\bea
\label{pr1}
T^{1/2}\left[\E(J_T^h)-  J\right] \to 0 \quad {\rm as}\quad T\to\infty.
\eea
\end{thm}
\begin{rem}
{\rm We call an estimator $J_T^h$ of $J$ {\it asymptotically unbiased of the order of\/}
$T^\be$,  $\be>0$ if $\lim_{T \to \f}T^\be[E (J_T^h) - J] = 0.$
Thus, Theorem \ref{AU3} states that the statistic $J_T^h$ is an asymptotically
unbiased estimator for $J$ of the order of $T^{1/2}$.}
\end{rem}

\sn{Consistency}

Recall that an estimator $J_T^h$ of $J$ is said to be
(a) consistent if $J_T^h\to J$ in probability as $T\to\f$,
(b) mean square consistent if $\E(J_T^h-J)^2 \to 0$ as $T\to\f$,
(c) $\sqrt{T}$-consistent in the mean square sense if
$\E\left([\sqrt{T}(J_T^h-J)]^2\right) =O(1)$ as $T\to\f$,
%Given $0 < \al < 1$, we define an estimator $J_T^h$ of $J(f)$ to be
%$T^\al$-consistent if $T^\al (J_T^h - J(f)) \to 0$ in probability as $T\to\f$.

To state the corresponding results we first introduce the following assumption.

\begin{asn}
\label{(A3)}
%{\bf Assumption (A3).}
{\rm The filter $a(\cdot)$ and the generating kernel $\widehat g(\cdot)$
are such that
\begin{equation}
\label{eq:CLT condition}
\nonumber
a(\cdot)\in L^p(\Lambda)\cap L^2(\Lambda),\quad \widehat g(\cdot)\in L^q(\Lambda)
\quad {\rm with}\quad 1\le p,q\le 2,\quad \frac{2}{p}+\frac{1}{q}\ge \frac{5}{2}.
\end{equation}}
\end{asn}

We begin with results on the asymptotic behavior of the variance
$\Var (J_T^h) = \E (J_T^h -\E (J_T^h))^2$.
The proof of the next theorem can be found in Ginovyan and Sahakyan \cite{GS2019}.
\begin{thm}
\label{CN1}
Let the functionals $J:=J(f,g)$ and $J_T^h: = J(I^h_{T},g)$ be defined
by (\ref {t1}) and (\ref{t7}), respectively. Then under
Assumptions \ref{(T)} and \ref{(A3)} the following asymptotic relation holds:
\bea
\label{cn2}
\lim_{T \to \f}T\Var (J_T^h) = \sigma^2_h(J),
\eea
where
\beq
\label{tsigma}
\sigma^2_h(J):=4\pi e(h)\int_{\Lambda} f^2(\lambda) g^2(\lambda) d\lambda
+\kappa_4 e(h)\left[\int_{\Lambda} f(\lambda) g(\lambda) d\lambda\right]^2.
\eeq
Here $\kappa_4$ is the fourth cumulant of $\xi(1)$, and
\beq
\label{eh}
e(h):=\frac{H_4}{H_2^2}= \int_0^1h^4(t)dt \left(\int_0^1h^2(t)dt\right)^{-2}.
\eeq
\end{thm}
From Theorems \ref{AU1}--\ref{CN1} we infer the following result.
\begin{thm}
\label{CN2} The following assertions hold.
\begin{itemize}
\item[(a)]
Under Assumptions \ref{(T)}, \ref{(A1)} and \ref{(A3)} the statistic $J_T^h$ is
a mean square consistent estimator for $J$.
\item[(b)]
Under Assumptions \ref{(T)}, \ref{(A2')} (or \ref{(A2)}) and \ref{(A3)} the statistic $J_T^h$ is
a $\sqrt{T}$-consistent in the mean square sense estimator for $J$.
\end{itemize}
\end{thm}

\sn{Asymptotic normality}

The next result contains sufficient conditions for functional $J_T^h$
to obey the central limit theorem (CLT), and was proved in Ginovyan and Sahakyan \cite{GS2019}.

\begin{thm}[CLT]
\label{T-CLT}
Let $J:=J(f,g)$ and $J_T^h: = J(I^h_{T},g)$ be defined by (\ref {t1}) and (\ref{t7}),
respectively. Then under Assumptions \ref{(T)} and \ref{(A3)} the functional $J_T^h$
obeys the central limit theorem. More precisely, we have
\bea
\label{pr2}
T^{1/2}\left[J_T^h - \E(J_T^h)\right] \ConvD \eta\quad {\rm as}\quad T\to\infty,
\eea
where the symbol $\ConvD$ stands for convergence in distribution, and
$\eta$ is a normally distributed random variable with mean zero and variance
$\si^2_h(J)$ given by \eqref{tsigma} and \eqref{eh}.
\end{thm}

Taking into account the equality
\bea
\label{pt1}
T^{1/2}\left[J_T^h - J\right]=T^{1/2}\left[\E(J_T^h)-  J\right]
+ T^{1/2}\left[J_T^h - \E(J_T^h)\right],
\eea
as an immediate consequence of Theorems \ref{AU3} and \ref{T-CLT},
we obtain the next result that contains sufficient conditions for a
simple 'plug-in' statistic $J(I^h_{T})$ to be an asymptotically normal
estimator for a linear spectral functional $J$. %given by (\ref{t1}).

\begin{thm}
\label{TT1}
Let the functionals $J:=J(f,g)$ and $J_T^h: = J(I^h_{T},g)$ be defined
by (\ref {t1}) and (\ref{t7}), respectively.
Then under Assumptions \ref{(T)}, \ref{(A2')} (or \ref{(A2)}) and \ref{(A3)} the statistic
$J_T^h$ is an asymptotically normal estimator for functional $J$.
More precisely, we have
\bea
\label{t8}
T^{1/2}\left[J_T^h-  J\right] \ConvD \eta\quad {\rm as}\quad T\to\infty,
\eea
where $\eta$ is as in Theorem \ref{T-CLT}, that is, $\eta$ is a normally
distributed random variable with mean zero and variance $\si^2_h(J)$
given by \eqref{tsigma} and \eqref{eh}.
\end{thm}

\begin{rem}
{\rm Notice that if the underlying process $X(u)$ is Gaussian,
then in formula (\ref{tsigma}) we have only the first term.
Using the results from Ginovyan \cite{G1994} and Ginovyan and
Sahakyan \cite{GS2005, GS2007}, it can be shown that in this case Theorem \ref{TT1}
is true under Assumptions \ref{(T)} and \ref{(A3)}.}
\end{rem}
\begin{exa}[Estimation of covariance function]
{\rm Assume that $X(t)$ is a c.t.\ process, and let
$g(\lambda) = e^{iu\lambda}$, then
$$J(f) =
\int_{-\f}^\f  e^{iu\lambda} f(\lambda)\,d\lambda : = r(u).$$
\noindent
Thus, in this special case our problem becomes to the estimation of the
the covariance function $r(u) = \mathbb{E}[X(t+u)X(t)]$ of the process $X(t)$.
By Theorem \ref{TT1} the simple "plug-in" statistic
$$
\widehat J_T^h = \widehat r_T(u) =
\int_{-\f}^\f e^{iu\lambda}I_T^h(\lambda)\,d\lambda
%= \frac1T\int_0^{T-|u|} X(t)X(t+u)dt
$$
\noindent
is asymptotically normal estimator for $r(u)$ with asymptotic variance
$$ \sigma_u^2 \, = \,
4\pi e(h)\int_{-\f}^\f \cos^2(u\la)f^2(\la)\,d\la,$$
where $e(h)$ is given by \eqref{eh}.}
\end{exa}

\begin{exa}[Estimation of spectral function]
{\rm Assume that $X(t)$ is a d.t.\ process, and let
$g(\lambda) = \chi_{[0,\mu]}(\lambda)$ be the indicator
of an interval $[0,\mu]$, then
$$J(f) =
\int_{-\pi}^\pi \chi_{[0,\mu]}(\lambda) f(\lambda)\,d\lambda
= \int_0^\mu f(\lambda)\,d\lambda : = F(\mu).$$

\noindent
Thus, in this case the estimand functional is the
spectral function $F(\mu)$ of the process $X(u),$
and by Theorem \ref{TT1} the simple "plug-in" statistic
\beq
\label{edf}
\widehat J_T^h = \widehat F_T(\mu) = \int_0^\mu I_T^h(\lambda)\,d\lambda
%=\frac1{2\pi T}\int_0^\mu\left|\sum_{t=1}^T X(t)e^{-i\lambda t}\right|^2
\eeq
\noindent
is asymptotically normal estimator for $F(\mu)$ with asymptotic variance
$$\sigma^2(\mu) =  4\pi e(h) \int_0^\mu f^2(s)\,ds,$$
where $e(h)$ is given by \eqref{eh}.}
\end{exa}

\s{Parametric Estimation Problem} %: The Whittle procedure}
\label{PW}
We assume here that the spectral density $f(\la)$ belongs to a given parametric family
of spectral densities $\mathcal{F}:=\{f(\lambda,\theta): \, \theta\in\Theta\}$,
where $\theta:=(\theta_1, \ldots, \theta_p)$ is an unknown parameter and
$\Theta$ is a subset in the Euclidean space $\mathbb{R}^p$.
The problem of interest is to estimate the unknown parameter $\theta$ on the basis
of the tapered data \eqref{t2}, %$\{h(t)X(t), \, 0\le t\le T\}$,
and investigate the asymptotic (as $T \to \infty$) properties of the suggested estimators,
%(consistency, asymptotic normality, etc.),
depending on the dependence (memory) structure of the model $X(t)$ and the
smoothness of its spectral density $f$.

There are different methods of estimation: maximum likelihood,
Whittle, minimum contrast, etc.
Here we focus on the Whittle method.
%, which is based on the smoothed periodogram analysis in the frequency domain.

\sn{The Whittle estimation procedure}

The Whittle estimation procedure, originally devised for d.t.\ short memory
stationary processes, is based on the smoothed periodogram analysis on a frequency
domain, involving approximation of the likelihood function and asymptotic properties
of empirical spectral functionals (see Whittle \cite{W}).
The Whittle estimation method since its discovery has played
a major role in the asymptotic theory of parametric estimation in the frequency domain,
and was the focus of interest of many statisticians.
Their aim was to weaken the conditions needed to guarantee the validity of the
Whittle approximation for d.t.\ short memory models,
to find analogues for long and intermediate memory models,
to find conditions under which the Whittle estimator is asymptotically equivalent
to the exact maximum likelihood estimator, and to extend the procedure to the
c.t.\ models and random fields.

For the d.t.\ case, it was shown that for Gaussian and linear
stationary models the Whittle approach leads to consistent and asymptotically
normal estimators under short, intermediate and long memory assumptions.
Moreover, it was shown that in the Gaussian case the Whittle
estimator is also asymptotically efficient in the sense of Fisher
(see, e. g., Dahlhaus \cite{D}, Dzhaparidze \cite{Dz1}, Fox and Taqqu \cite{FT1},
Giraitis and Surgailis \cite{GSu},
Guyon \cite{Gu}, Heyde and Gay \cite{HG}, Taniguchi and Kakizawa \cite{TK},
Walker \cite{Wal}, and references therein).

%Continuous versions of Whittle estimation procedure have been
For c.t.\ models, the Whittle estimation procedure has been
considered, for example, in Anh et al. \cite{ALS2}, Avram et al. \cite{AvLS},
Casas and Gao \cite{CG}, Dzhaparidze \cite{Dz1},
Dzhaparidze and Yaglom \cite{DY}, Gao \cite{Go},
Gao et al. \cite{GAHT}, Leonenko and Sakhno \cite{LS},
Tsai and Chan \cite{TC},
where can also be found additional references. In this case, it was proved that
the Whittle estimator is consistent and asymptotically normal.

The Whittle estimation procedure based on the d.t.\ tapered data
has been studied in Alomari et al. \cite{ALRST}, Dahlhaus \cite{D1},
Dahlhaus and K\"unsch \cite{DK}, Guyon \cite{Gu}, Lude\~na and Lavielle \cite{LL}.
In the case where the underlying model is a L\'evy-driven c.t.\
linear process with possibly unbounded or vanishing spectral density
function, consistency and asymptotic normality of the Whittle estimator
was established in Ginovyan \cite{G2020e}.

To explain the idea behind the Whittle estimation procedure, assume for simplicity
that the underlying process $X(t)$ is a d.t.\ Gaussian process,
and we want to estimate the parameter $\theta$ based on the sample
$X_T:=\{X(t), \, t=1,\ldots, T\}$. A natural approach is to find the maximum likelihood
estimator (MLE) $\widehat\theta_{T,MLE}$ of $\theta$, that is, to maximize the
likelihood function, or to minimize the $-1/T\times$log-likelihood function $L_T(\theta)$,
which in this case takes the form:
$$L_T(\theta):=\frac12\ln2\pi+\frac1{2T}\ln\det B_T(f_\theta)+\frac1{2T}X'_T[B_T(f_\theta)]^{-1}X_T,$$
where $B_T(f_\theta)$  is the Toeplitz matrix generated by $f_\theta$.
Unfortunately, the above function is difficult to handle, and no explicit expression
for the estimator $\widehat\theta_{T,MLE}$ is known (even in the case of simple models).
An approach, suggested by P. Whittle, called the Whittle estimation procedure,
is to approximate the term $\ln\det B_T(f_\theta)$
by $\frac T2\int_{-\pi}^\pi\ln f_\theta(\la)d\la$ and the inverse matrix $[B_T(f_\theta)]^{-1}$
by the Toeplitz matrix $B_T(1/f_\theta)$.
This leads to the following approximation of the log-likelihood function $L_T(\theta)$,
introduced by Whittle \cite{W}, and called Whittle functional:
$$L_{T,W}(\theta)=\frac 1{4\pi}\int_{-\pi}^\pi\left[\ln f_\theta(\la)
+\frac{I_{T}(\la)}{f_\theta(\la)}\right]\, d\la,$$
where $I_{T}(\lambda)$ is the ordinary periodogram of the process $X(t)$.

Now minimizing the Whittle functional $L_{T,W}(\theta)$ with respect to
$\theta$, we get the Whittle estimator $\widehat\theta_{T}$ for $\theta$.
It can be shown that if
$$ T^{1/2}(L_T(\theta)- L_{T,W}(\theta)\to0 \q {\rm as}\q n\to\f \q {\rm in \,\, probability,}
$$
then the MLE $\widehat\theta_{T,MLE}$ and the Whittle estimator $\widehat\theta_{T}$
are asymptotically equivalent in the sense that $\widehat\theta_{T}$ also is
consistent, asymptotically normal and asymptotically Fisher-efficient
(see, e.g., Coursol and Dacunha-Castelle \cite{DC},
Dzhaparidze \cite{Dz1}, and Dzhaparidze and Yaglom \cite{DY}).
%(see, e.g., Dzhaparidze \cite{Dz1}, Dzhaparidze and Yaglom \cite{DY}).

In the continuous context, the Whittle procedure of estimation of a spectral
parameter $\theta$ based on the sample $X_T:=\{X(t), \, 0\leq t\leq T\}$
is to choose the estimator $\widehat\theta_{T}$ to minimize the weighted
Whittle functional:
\beq
\label{pe9a}
U_{T}(\theta): =\frac1{4\pi}\int_{\mathbb{R}}\left[\ln f(\la, \theta) +
\frac{I_{T}(\la)}{f(\la, \theta)}\right]\cd w(\la) \, d\la,
\eeq
where $I_{T}(\la)$ is the continuous periodogram of $X(t)$,
and $w(\la)$ is a weight function ($w(-\la)=w(\la)$, $w(\la)\ge0$,
$w(\la)\in L^1(\mathbb{R})$) for which the integral in (\ref{pe9a})
is well defined.
An example of common used weight function is $w(\la)=1/(1+\la^2)$.

The Whittle procedure of estimation of a spectral parameter $\theta$
based on the tapered sample \eqref{t2} is to choose the estimator
$\widehat\theta_{T,h}$ to minimize the weighted tapered Whittle functional:
\beq
\label{pe9}
U_{T,h}(\theta): =\frac1{4\pi}\int_{\Lambda}\left[\log f(\la, \theta) +
\frac{I^h_{T}(\la)}{f(\la, \theta)}\right]\cd w(\la) \, d\la,
\eeq
where $I^h_{T}(\la)$ is the tapered periodogram of $X(t)$, given by \eqref{t6},
and $w(\la)$ is a weight function for which the integral in (\ref{pe9})
is well defined. Thus,
\beq
\label{pee9}
\widehat\theta_{T,h}: =\underset{\theta\in\Theta}{\rm Arg \, min} \, U_{T,h}(\theta).
\eeq

\sn{Asymptotic properties of the Whittle estimator}

To state results involving properties of the Whittle estimator,
we first introduce the following set of assumptions.

\begin{asn}
\label{(B1)}
%{\bf Assumption (B1).}
\rm{The true value $\theta_0$ of the parameter $\theta$ belongs to a compact set
$\Theta$, which is contained in an open set $S$ in the $p$-dimensional
Euclidean space $\mathbb{R}^p$,
%the true value of the parameter $\theta_0$ is in the interior of $\Theta$,
and $f(\la,\theta_1)\neq f(\la,\theta_2)$ whenever $\theta_1\neq \theta_2$
almost everywhere in $\Lambda$ with respect to the Lebesgue measure.}
\end{asn}

\begin{asn}
\label{(B2)}
%{\bf Assumption (B2).}
{\rm The functions $f(\la,\theta)$, $f^{-1}(\la,\theta)$ and
$\frac{\partial}{\partial\theta_k}f^{-1}(\la,\theta)$, $k=1,\ldots,p$,
are continuous in $(\la,\theta)$.}
\end{asn}

\begin{asn}
\label{(B3)}
%{\bf Assumption (B3).}
{\rm The functions $f:=f(\la,\theta)$ and
$g:=w(\la)\frac{\partial}{\partial\theta_k}f^{-1}(\la,\theta)$
satisfy Assumptions \ref{(A2)} or \ref{(A3)} %(see Section \ref{NP})
for all $k=1,\ldots,p$ and $\theta\in\Theta$.}
\end{asn}

\begin{asn}
\label{(B4)}
%{\bf Assumption (B4).}
{\rm The functions $a:=a(\la,\theta)$ and $b:=\widehat g$, where $g$ is as in
Assumption \ref{(B3)}, satisfy Assumption \ref{(A1)}.} % (see Section \ref{NP}).}
\end{asn}

\begin{asn}
\label{(B5)}
%{\bf Assumption (B5).}
{\rm The functions
$\frac{\partial^2}{\partial\theta_k\partial\theta_j}f^{-1}(\la,\theta)$ and
$\frac{\partial^3}{\partial\theta_k\partial\theta_j\partial\theta_j}f^{-1}(\la,\theta)$,
$k,j,l=1,\ldots, p$, are continuous in $(\la,\theta)$ for $\la\in\Lambda$,
$\theta\in N_\de(\theta_0)$, where $N_\de(\theta_0):=\{\theta: \, |\theta-\theta_0|<\de\}$
is some neighborhood of $\theta_0$.}
\end{asn}

\begin{asn}
\label{(B6)}
%{\bf Assumption (B6).}
{\rm The matrices
\bea
\label{W03}
W(\theta):=\|w_{ij}(\theta)\|_{i,j=1,\ldots,p}, \q
A(\theta):=\|a_{ij}(\theta)\|_{i,j=1,\ldots,p},\q
B(\theta):=\|b_{ij}(\theta)\|_{i,j=1,\ldots,p}
\eea
are positive definite, where
\bea
\label{W3}
w_{ij}(\theta)&=&  \frac1{4\pi}\int_{\Lambda}
\frac{\partial}{\partial\theta_i}\ln f(\la, \theta)
\frac{\partial}{\partial\theta_j}\ln f(\la, \theta)w(\la)d\la,\\
\label{W4}
a_{ij}(\theta)&=& \frac1{4\pi}\int_{\Lambda}
\frac{\partial}{\partial\theta_i}\ln f(\la, \theta)
\frac{\partial}{\partial\theta_j}\ln f(\la, \theta)w^2(\la)d\la,\\
\label{W5}
b_{ij}(\theta)&=& \frac{\kappa_4}{16\pi^2}\int_{\Lambda}
\frac{\partial}{\partial\theta_i}\ln f(\la, \theta)w(\la)d\la
\int_{\mathbb{R}} \frac{\partial}{\partial\theta_j}\ln f(\la, \theta)w(\la)d\la,
\eea
and $\kappa_4$ is the fourth cumulant of $\xi(1)$.}
\end{asn}

%\n{\it Consistency of the Whittle estimator.}
The next theorem contains sufficient
conditions for Whittle estimator to be consistent (see Ginovyan \cite{G2020e}).
\begin{thm}
\label{CWT}
Let $\widehat\theta_{T,h}$ be the Whittle estimator defined by \eqref{pee9} and let
$\theta_0$ be the true value of parameter $\theta$. Then, under Assumptions
\ref{(B1)}--\ref{(B4)} and \ref{(T)}, % (see Section \ref{H-6.1}),
the statistic $\widehat\theta_{T,h}$ is a consistent estimator
for $\theta$, that is, $\widehat\theta_{T,h}\to \theta_0$ in probability as $T\to\f$.
\end{thm}

\n
%{\it Asymptotic normality of the Whittle estimator.}
Having established the consistency of the Whittle estimator $\widehat\theta_{T,h}$,
we can go on to obtain the limiting distribution of $T^{1/2}\left(\widehat\theta_{T,h}-\theta_0\right)$
in the usual way by applying the Taylor's formula, the mean value theorem, and Slutsky's arguments.
Specifically we have the following result, showing that under the above assumptions,
the Whittle estimator $\widehat\theta_{T,h}$ is asymptotically normal
(see Ginovyan \cite{G2020e}).

\begin{thm}
\label{TAN}
Suppose that Assumptions \ref{(B1)}--\ref{(B6)} and \ref{(T)} %, (see Section \ref{H-6.1})
are satisfied.
%\ref{T1} with $f=f(\la; \theta)$ and $g=w(\la)\frac{\partial}{\partial\theta_i}f^{-1}(\la, \theta)$ $(i=1,\ldots,p)$
Then the Whittle estimator $\widehat\theta_{T,h}$ of an unknown spectral parameter
$\theta$ based on the tapered data \eqref{t2} is asymptotically normal. More precisely,
we have
\bea
\label{aW1}
T^{1/2}\left(\widehat\theta_{T,h}-\theta_0\right)\ConvD N_p\left(0, e(h)\G(\theta_0)\right)
\quad {\rm as}\quad T\to\infty,
\eea
where $N_p(\cdot,\cdot)$ denotes the $p$-dimensional normal law, \, $\ConvD$
stands for convergence in distribution,
\bea
\label{aW2}
\G(\theta_0) = W^{-1}(\theta_0)\left(A(\theta_0)+B(\theta_0)\right)W^{-1}(\theta_0),
\eea
where the matrices $W$, $A$ and $B$ are defined in (\ref{W03})-(\ref{W5}),
and the tapering factor $e(h)$ is given by formula \eqref{eh}.
\end{thm}

\begin{rem}
{\rm
In the d.t.\ case as a weight function we take $w(\la)\equiv 1$,
and the matrices $A(\theta_0)$ and $W(\theta_0)$ coincide (see (\ref{W03}) -- (\ref{W4})).
So, in this case, formula \eqref{aW2} becomes
$\G(\theta_0) = W^{-1}(\theta_0)\left(W(\theta_0)+B(\theta_0)\right)W^{-1}(\theta_0)$.
If, in addition, the underlying process is Gaussian  ($\kappa_4=0$, and hence $B(\theta_0)=0$), and the taper $h$ is chosen
so that the tapering factor $e(h)$ is equal to one, then we have
$\G(\theta_0) = W^{-1}(\theta_0)$, that is, the Whittle estimator $\widehat\theta_{T,h}$
is Fisher-efficient.}
\end{rem}

\section{Goodness-of-fit tests}
\label{GF}

In this section we consider the following problem of hypotheses testing.

Based on the tapered sample $\mathbf{X}_T^h$ given by \eqref{t2},
we want to construct goodness-of-fit tests for
testing a hypothesis $H_0$ that the spectral density of the process $X(t)$
has the specified form $f(\la)$. We will distinguish the following two cases.
\begin{itemize}
\item[a)]
The hypothesis $H_0$ is simple, that is, the hypothetical spectral density
$f(\la)$ of $X(t)$ does not depend on unknown parameters.
\item[b)]
The hypothesis $H_0$ is composite, that is, the hypothetical spectral density
$f(\la)$ of $X(t)$ depends on an unknown $p$--dimensional vector parameter
${\theta} =(\theta_1, \ldots, \theta_p)$,
that is, $f(\la)=f(\la, {\theta})$, $\la\in\Lambda$,
${\theta}\in \Theta\subset\mathbb{R}^p$.
%, where $S$ is an open set of the Euclidean space $\Lambda^p$.
\end{itemize}

The above stated problem in the non-tapered case has been considered by many authors
for different models.
For instance, for independent observations, the problem was considered in
Chernov and Lehman \cite{CL}, Chibisov \cite{Ch}, Cramer \cite{C}, and
Dzhaparidze and Nikulin \cite{DN}. %and Kendall and Stuart \cite{KS}.
For observations generated by d.t.\ Gaussian stationary processes
it was considered in Dzhaparidze \cite{Dz}, Ginovyan \cite{G2003g},
Hannan \cite{H}, and Osidze \cite{O1,O2}.
For c.t.\ Gaussian stationary observations, the problem
was discussed in Ginovyan \cite{G2018g} and Osidze \cite{O1,O2}.
For tapered case the problem has been considered in Ginovyan \cite{G2020g}.

To test the hypothesis $H_0$, similar to the non-tapered case, it is natural to
introduce a measure of divergence (disparity) of the hypothetical and empirical
spectral densities, and construct a goodness-of-fit test based on the distribution
of the chosen measure (see, e.g., Dzhaparidze \cite{Dz},
Ginovyan \cite{G2003g,G2018g,G2020g}, Hannan \cite{H}, and Osidze \cite{O1,O2}).

\subsection{A Goodness-of-fit test for simple hypothesis}

We first consider the relatively easy case a) of a simple hypothesis $H_0$.
As a measure of divergence of the hypothetical spectral density
$f(\la)$ and the tapered empirical spectral density $I^h_T(\lambda)$,
we consider the $m$--dimensional random vector
\beq \label{1-0}
{\Phi}^h_{T}:= \bigl(\Phi^h_{1T},\ldots,\Phi^h_{mT}\bigr)
\eeq
with elements
\beq \label{1-1}
\Phi^h_{jT}:  \, = \,\Phi_{jT}(\mathbf{X}_T^h)  \, = \,
\frac{\sqrt{T}}{\sqrt{4\pi e(h)}}\int_{\Lambda}
\biggl[\frac{I^h_T(\la)}{f(\la)}-1\biggr]\,\I_j(\la)d\la,
\quad j=1,2,\ldots, m,
\eeq
where $e(h)$ is as in \eqref{eh} and $\{\I_j(\la), \, j=1,2,\ldots, m\}$
is some orthonormal system on $\Lambda$:
\beq
\label{o1}
%\frac1{2\pi}
\int_{\Lambda}\I_k(\la)\I_j(\la)\,d\la =\de_{kj}
=\left \{
\begin{array}{ll}
1& \mbox{for $k=j$},\\
\\[-4mm]
0& \mbox{for $k\neq j$}.
\end{array}
\right.
\eeq
%In Lemma 3.2 of Ginovyan \cite{G2020g}
In Ginovyan \cite{G2020g} it was shown that under wide conditions on
$f(\la)$ and $\I_j(\la)$, the random vector
${\Phi}^h_{T}$ in \eqref{1-0} -- \eqref{1-1} has asymptotically $N(0, I_m)$--normal
distribution as $T\to\f$, where $I_m$ is the $m\times m$ identity matrix.
%, and the components $\Phi^h_{kT}$ and $\Phi^h_{jT}$
%of ${\Phi}^h_{T}$ are asymptotically uncorrelated for $k\ne j$.
Therefore in the case of simple hypothesis $H_0$,  we can use the statistic
\beq
\label{1-3}
S_{T}^h=S_{T}(\mathbf{X}_T^h):= {\Phi}'_{T}(\mathbf{X}_T^h){\Phi}_{T}(\mathbf{X}_T^h)
=\sum_{j=1}^m \Phi^2_{jT}(\mathbf{X}_T^h),
\eeq
which for $T\to\f$ will have a $\chi^2$--distribution with $m$ degrees of freedom.
% (see \cite{Dz-1}, \cite{H}).

Thus, fixing an asymptotic level of significance $\al$ we can
consider the class of goodness-of-fit tests for testing the simple
hypothesis $H_0$ about the form of the spectral density $f$ with
asymptotic level of significance $\al$ determined by critical
regions of the form:
\beq
\label{13}
\nonumber
\{\mathbf{X}_T^h:\, \, S_{T}(\mathbf{X}_T^h)> d_\al\},
\eeq
where $S_{T}(\mathbf{X}_T^h)$ is given by (\ref{1-3}), and $d_\al$ is the $\al$-quantile
of $\chi^2$--distribution with $m$ degrees of freedom, that is,  $d_\al$ is
determined from the condition:
\beq
\label{14}
\nonumber
P(\chi^2> d_\al) = \int_{d_\al}^\f k_m(x)\, dx =\al,
\eeq
where $k_m(x)$ is the density of $\chi^2$--distribution with $m$ degrees of freedom.

The next theorem contains sufficient conditions for statistic $S_{T}^h$,
given by (\ref{1-3}), to have a limiting (as $T\to\f$) $\chi^2$--distribution
with $m$ degrees of freedom (see Ginovyan \cite{G2020g}).
\begin{thm}
\label{ST1}
Let the spectral density $f(\la)$ and the orthonormal functions
$\{\I_j(\la), \, j=1,2,\ldots, m\}$ be such that $(f, g_j)\in (\mathcal{H})$
for all $j=1,2,\ldots, m$, where $g_j=\I_j/f$ (see Assumption \ref{(A2')}).
Then under Assumption \ref{(T)}
the limiting (as $T\to\f$) distribution of the statistic
$S_{T}^h=S_{T}(\mathbf{X}_T^h)$ given by (\ref{1-3})
is a $\chi^2$--distribution with $m$ degrees of freedom.
\end{thm}
\begin{rem}
\label{r5.1}
{\rm For the non-tapered case, for observations generated by  d.t.\
short-memory Gaussian stationary processes the result of Theorem
\ref{ST1} was first proved in  Hannan \cite{H} (p. 94) (see, also,
Dzhaparidze \cite{Dz} and Osidze \cite{O1,O2}).
In the case where the spectral density has singularities (zeros and/or poles),
the result for d.t.\ processes was proved in Ginovyan \cite{G2003g}.
The non-tapered counterpart of Theorem \ref{ST1} for c.t.\ processes
was proved in Ginovyan \cite{G2018g}.}
\end{rem}

\subsection{A Goodness-of-fit test for composite hypothesis}

Now we consider the case of composite hypothesis $H_0$, and assume that
the hypothetical spectral density $f=f(\lambda,\theta)$ is known with the
exception of a vector parameter
$\theta:=(\theta_1, \ldots, \theta_p)\in\Theta\subset\mathbb{R}^p$.
In this case, the problem of construction of goodness-of-fit tests becomes
more complex, because we first have to choose
an appropriate statistical estimator $\widehat{\ta}_T$ for the unknown
parameter $\theta$, constructed on the basis of the tapered sample \eqref{t2}.
It is important to remark that in this case the limiting distribution of the test
statistic will change in accordance with properties of an estimator of ${\theta}$,
and generally will not be a $\chi^2$--distribution.

For testing a composite hypothesis $H_0$, we again can use a statistic of type (\ref{1-3}),
but with a statistical estimator $\widehat{{\ta}}_T$ instead of unknown ${\ta}$.
The corresponding statistic can be written as follows:
\beq \label{1-4}
S_{T}^h(\widehat{\ta}_T)= S_{T}(\mathbf{X}_T^h, \widehat{\ta}_T):
= {\Phi}'_{T}(\mathbf{X}_T^h, \widehat{\ta}_T)
{\Phi}_{T}(\mathbf{X}_T^h, \widehat{\ta}_T) =\sum_{j=1}^m \Phi^2_{jT}(\mathbf{X}_T^h, \widehat{\ta}_T),
\eeq
where now
\beq \label{1-0c}
\nonumber
{\Phi}^h_{T}(\mathbf{X}_T^h, \widehat{\ta}_T):= \bigl(\Phi_{1T}(\mathbf{X}_T^h, \widehat{\ta}_T),
\ldots, \Phi_{mT}(\mathbf{X}_T^h, \widehat{\ta}_T)\bigr)
\eeq
with elements
\beq \label{1-1c}
\Phi_{jT}(\mathbf{X}_T^h, \widehat{\ta}_T):  \, = \,
\frac{\sqrt{T}}{\sqrt{4\pi e(h)}}\int_{\Lambda} \biggl[\frac{I^h_T(\la)}{f(\la,\widehat{\ta}_T)}-1\biggr]
\,\I_j(\la)d\la, \quad j=1,2,\ldots, m.
\eeq

So, we must choose an appropriate statistical estimator $\widehat{\ta}_T$ for
unknown ${\ta}$, and determine the limiting distribution of the statistic (\ref{1-4}).
Then, having the limiting distribution of the statistic (\ref{1-4}), for given level of
significance $\al$ we can consider the class of goodness-of-fit tests for testing
the composite hypothesis $H_0$ about the form of the spectral density $f$ with asymptotic
level of significance $\al$ determined by critical regions of the form:
\beq
\label{15}
\nonumber
\{\mathbf{X}_T^h:\, \, S_{T}(\mathbf{X}_T^h, \widehat{\ta}_T)> d_\al\},
\eeq
where $d_\al$ is the $\al$-quantile of the limiting distribution
of the statistic (\ref{1-4}), that is,  $d_\al$ is determined
from the condition:
\beq
\label{16}
\nonumber
 \int_{d_\al}^\f \widehat k_m(x)\, dx =\al,
\eeq
where $\widehat k_m(x)$ is the density of the limiting distribution of
$S^h_{T}(\widehat{\ta}_T)$ defined by (\ref{1-4}).

To state the corresponding result we first introduce the following set of assumptions:

%\begin{asn}
%\label{(GA1)}
%{\bf Assumption (A1).}
%{\rm The true value $\ta_0$ of the parameter  $\ta$ belongs to a
%bounded closed set $\Theta$ %contained in an open set $S$
%in the $p$-dimensional Euclidean space $\Lambda^p$.}
%\end{asn}

%\begin{asn}
%\label{(GA2)}
%{\bf Assumption (A2).}
%{\rm If $\ta_1$ and $\ta_2$ are two distinct points of $\Theta$,
%then $f(\la, \ta_1)\neq f(\la, \ta_2)$ almost everywhere in $\Lambda$
%with respect to the Lebesgue measure.}
%\end{asn}

\begin{asn}
\label{(GA3)}
%{\bf Assumption (A3).}
{\rm For $\ta\in\Theta$, $(f, g_j)\in (\mathcal{H})$
for all $j=1,2,\ldots, m$, where $f:=f(\la, \ta)$ and
$g_j:=\I_j(\la)/f(\la, \ta)$.}
\end{asn}

\begin{asn}
\label{(GA4)}
%{\bf Assumption (A4).}
{\rm For $\ta\in\Theta$, $(f, h_{kj})\in (\mathcal{H})$
for all $k=1,2,\ldots, p$ and $j=1,2,\ldots, m$, where $f:=f(\la, \ta)$ and
$h_{kj}:=\dfrac{\I_j(\la)}{f(\la, \ta)}\dfrac{\p}{\p\ta_k}\ln f(\la,\ta)$.}
\end{asn}

\begin{asn}
\label{(GA5)}
%{\bf Assumption (A5).}
{\rm The $(p\times p)$--matrix $\G(\ta_0)= ||\g_{kj}(\ta_0)||_{k, j=\ol{1,p}} \, $
with elements
\beq \label{2-1}
\g_{kj}(\ta_0):=
\frac{1}{4\pi}\int_{\Lambda}
\left[\frac{\p}{\p\ta_k}\ln f(\la,\ta)\right]_{\ta=\ta_0}
\left[\frac{\p}{\p\ta_j}\ln f(\la,\ta)\right]_{\ta=\ta_0}\,d\la
\eeq
is non-singular.}
\end{asn}

\begin{asn}
\label{(GA6)}
%{\bf Assumption (A6).}
{\rm There exists a $\sqrt T$--consistent estimator $\best$ for the
parameter $\ta$ such that the following asymptotic relation holds:
\beq \label{3-4}
\sqrt{T}\,(\best -\ta_0) - \G^{-1}(\ta_0)\De_T(\ta_0)= o_P(1),
\eeq
where $\G^{-1}(\ta_0)$ is the inverse of the matrix $\G(\ta_0)$ defined in Assumption \ref{(GA5)}, and
\beq
\label{3-3a}
\De_T(\ta):=\De^h_T(\ta)=\bigl(\De_{1T}(\ta),\ldots, \De_{pT}(\ta)\bigr)
\eeq
is a $p$-dimensional random vector with components
\beq
\label{3-3}
\De_{kT} (\ta):  \, =
\, \frac{\sqrt{T}}{\sqrt{4\pi e(h)}}\int\limits_{\Lambda}
\biggl[\frac{I^h_T(\la)}{f(\la,\ta)}-1\biggr]\,
\frac{\partial}{\partial\theta_{k}}\ln f(\la,\ta)\,d\la,
\quad k=1, \ldots, p.
\eeq
The term $o_P(1)$ in \eqref{3-4} tends to zero in probability as $T\to\f$.
(Recall that an estimator $\best$ for $\ta$ is called $\sqrt T$--consistent
if $\sqrt T(\best-\ta)$ is bounded in probability).}
\end{asn}
\begin{rem}
\label{RH2}
{\rm As an estimator $\best$ for $\ta$ satisfying \eqref{3-4}
can be considered minimum contrast estimators (in particular, the Whittle estimator)
based on the tapered data.
Minimum contrast estimators based on the tapered data for
d.t.\ processes have been studied in Dahlhaus \cite{D1,D3,D2},
for Gaussian c.t.\ processes in Ginovyan \cite{G2020e},
and for some classes of c.t.\ non-Gaussian processes in Alomari et al \cite {ALRST}.}
\end{rem}

Let $B(\ta):= B^h(\ta)=||b_{jk}(\ta)||_{j=\ol{1,m}, \, k=\ol{1,p}}, \,$
be a $\, (m\times p)$--matrix with elements
\beq \label{2-4}
b_{jk}(\ta):= \frac{1}{\sqrt{4\pi e(h)}}\int_{\Lambda} \I_j(\la)
\frac{\p}{\p\ta_k}\ln f(\la,\ta)\,d\la,
\eeq
where $\I_j(\la)$ $(j=1,2,\ldots, m)$ are the functions from
(\ref{1-1}), and $e(h)$ is as in \eqref{eh}.

The following theorem was proved in Ginovyan \cite{G2020g}.
\begin{thm}\label{CT1}
Under Assumptions \ref{(T)}, \ref{(B1)} and \ref{(GA3)}--\ref{(GA6)}  the limiting distribution
(as $T\to\f$) of the statistic $S_{T}(\mathbf{X}_T^h, \widehat{\ta}_T)$
given by (\ref{1-4}), coincides with the distribution of the random variable
\beq \label{2-6}
\nonumber
\sum_{j=1}^{m-p} \xi^2_{j}+
\sum_{j=1}^p \nu_{j} \, \xi^2_{m-p+j},
\eeq
where $\xi_j$, $j=1,2,\ldots, m,$ are iid $N(0,1)$ random variables,
while the numbers $\nu_k$ $(0\le\nu_k<1)$, $k=1,2, \ldots, p$, are the
roots relative to $\nu$ of the following equation:
\beq \label{2-7}
\det\left[(1-\nu)\Gamma({\theta}_0)-
B'({\theta}_0)\,B({\theta}_0)\right]=0.
\eeq
\end{thm}
\begin{rem}
\label{r2}
{\rm For the non-tapered case, for independent observations the result of Theorem \ref{CT1}
was first obtained by Chernov and Lehman \cite{CL} (see, also, Chibisov \cite{Ch}).
For observations generated by d.t.\ short-memory Gaussian stationary
processes the result was stated in Osidze \cite{O1} (see, also,
Dzhaparidze \cite{Dz}). In the case where the spectral density has singularities,
the result for d.t.\ processes was proved in Ginovyan \cite{G2003g}.
The non-tapered counterpart of Theorem \ref{CT1} for c.t.\ processes
was proved in Ginovyan \cite{G2018g}.}
\end{rem}
%\begin{rem}
%\label{rn2}
%{\rm Comparing Theorems \ref{ST1} and \ref{CT1} with the corresponding
%results obtained in Ginovyan \cite{G7,G8} for the non-tapered case,
%we conclude that the limiting distributions of the statistic $S_{T}^h=S_{T}(\mathbf{X}_T^h)$ and $S_{T}(\mathbf{X}_T^h, \widehat{\ta}_T)$
%given by (\ref{1-3}) and (\ref{1-4}), respectively, do not depend
%on the taper factor .}
%\end{rem}
\begin{exa}
{\rm Let $X(t)$ be a d.t.\ Autoregressive Process of order $p$ (AR$(p)$),
that is, $X(t)$ is a stationary process with the spectral density
$
f(\la) = \frac{1}{2\pi}\cd \bigl|\al_p(e^{-i\la})\bigr|^{-2},
$
where $\al_p(z)=1-\ta_1z-\cdots-\ta_pz^p$. Consider the functions
$\I_j(\la)$ ($j=1,\ldots,m$) defined by
$\I_j(\la):=ce^{ij\la}\al_p(e^{-i\la})\bigl|\al_p(e^{-i\la})\bigr|^{-1}$
for $j>p$ and $\I_j(\la):=0$ for $j\leq p$, where $c$ is a normalizing
constant. As an estimator of $\ta:=(\ta_1, \ldots,\ta_p)$ consider the
Whittle estimator, and as a taper
the Tukey-Hanning taper function  $h(t)$ (see Remark \ref{TH}).
Then it is easy to check that the conditions of Theorem \ref{CT1} are
satisfied and $B(\ta)=0$. Therefore, the limiting (as $T\to\f$) distribution
of the statistic in (\ref{1-4})
%$S_{T}^h(\widehat{\ta}_T)$ given by (\ref{1-4})
is a $\chi^2$--distribution with $m-p$ degrees of freedom.}
\end{exa}

\section{Robustness to small trends of estimation}
\label{Rob}
In time series analysis, much of statistical inferences about unknown spectral
parameters or spectral functionals are concerned with the stationary
models, in which case it is assumed that the models are centered, or have constant
means.
In this section, we are concerned with the robustness of inferences, carried out on
a stationary models, possibly exhibiting long memory, contaminated by a small trend.
Specifically, let $\{X(t),\, t\in \mathbb{U}\}$ be a centered stationary
process possessing a spectral density $f_X(\lambda)$, $\lambda\in \Lambda$.
Assuming that either $f_X$ is known with the exception of a vector parameter
$\theta\in\Theta\subset\mathbb{R}^p$,
or $f_X$ is completely unknown and belongs to a given class $\mathcal{F}$,
we want to make inferences about $\theta$ or the value $J(f_X)$ of a given
functional $J(\cdot)$ at an unknown point $f_X\in\mathcal{F}$
in the case where the actual observed data are in the contaminated form:
\beq
\label{ri1}
Y(t)=X(t) + M(t), \quad 0\le t\le T,
\eeq
where $M(t)$ is a deterministic trend.

The process $X(t)$ is what we believe is being observed but in reality the
data are in the contaminated form $Y(t)$. %Assuming that the trend $M(t)$ is small,
%(that is, assuming that the major trend is removed from the model and a certain
%component that remains in the model has only minor effect),
In this case standard inferences can be carried on the basis of the
stationary model $X(t)$, and we are interested in question whether the
conclusions are robust against this kind of departure from the stationarity.

This problem for d.t.\ models was considered in Heyde and Dai \cite{HD}
(see also Taniguchi and Kakizawa \cite{TK}, Theorems 6.4.1 and 6.4.2).
For c.t.\ models it was studied in Ginovyan and Sahakyan \cite{GS2016}.

The results stated below show that if the trend $M(t)$ is 'small', then the
asymptotic properties of estimators of the parameter $\theta$ and the functional
$J(f)$, stated in Sections \ref{NP} and \ref{PW} for a stationary model $X(t)$,
remain valid for the contaminated model $Y(t)$, that is, both the parametric and
nonparametric estimating procedures are robust against replacing the stationary
model $X(t)$ by the non-stationary $Y(t)$. To this end, we first establish
an asymptotic relation between stationary and contaminated periodograms.
%Note that this result is true for a general stationary model $Y(t)$, and does
%not require from $Y(t)$ to have the specific form (\ref{lp}).

\subsection{A relation between stationary and contaminated periodograms}

The next result shows that a small trend of the form $|M(t)|\leq C|t|^{-\beta}$,
$\be>1/4$, does not effect the asymptotic properties of the empirical spectral
linear functionals of a periodogram.
Note that this result is of general nature,
and do not require from the model $X(t)$ to be linear. % have the specific form (\ref{lp}).

\begin{thm}
\label{T4}
Let $\{X(t), \, t\in \mathbb{U}\}$ be a stationary mean zero process,
$\{M(t), \, t\in \mathbb{U}\}$ be a deterministic trend, $Y(t)=X(t) + M(t)$,
and let $I^h_{TX}(\la)$ and $I^h_{TY}(\la)$ be the periodograms of $X(t)$
and $Y(t)$, respectively. Let $g(\la),$ $\la\in \Lambda$ be an even
integrable function.
If the trend $M(t)$ and the Fourier transform $a(t):=\widehat g(t)$
of $g(\la )$ are such that $M(t)$ is locally integrable on ${\mathbb R}$ and
\beq\label{re151}
|M(t)|\leq C|t|^{-\beta},\qquad |a(t)|\leq C|t|^{-\g}, \quad
t\in\Lambda,\quad 2\beta+\g>\frac 32,
\eeq
with some constants $C > 0$, $\g>0$ and  $\be>1/4$, then
\bea
\label{re9}
T^{1/2}\int_{\Lambda} g(\la)\left[I^h_{TY}(\la)-I^h_{TX}(\la) \right]d\la
\ConvP 0 \quad {\rm as}\quad T\to\infty,
\eea
where $\ConvP$ stands for convergence in probability,
provided that one of the following conditions holds:
\begin{itemize}
\item[(i)]
the process $X(t)$ has short or intermediate memory, that is,
the covariance function $r(t):=r_X(t)$ of $X(t)$ satisfies
$\ r\in L^1(\Lambda)$, and  $\be+\g>1$,
\item[(ii)]
the process $X(t)$ has long memory with %locally integrable
covariance function $r(t)$ satisfying
\beq\label{re99}
|r(t)| \le C|t|^{-\al}, \quad  t\in\Lambda,\quad \al+\g\ge\frac32
\eeq
with some constants  $C>0$, $0<\al\leq 1$,  and $\al+2\be>1$ if $\be<1<\g$.
\end{itemize}
\end{thm}

\begin{rem}
{\rm It is easy to check that the statement of Theorem \ref{T4} holds,
in particular, if the parameters $\al$, $\be$ and $\g$ satisfy the following
conditions:

 in the case  {\it (i)}: $\be>1/2,\ \g\geq1/2$,

 in the case {\it (ii)}: $\al\geq3/4, \ \be>3/8,\ \g\geq3/4$.}
\end{rem}
\begin{rem}
{\rm In the d.t.\ case, Theorem \ref{T4}
(with additional conditions $\g=1$ in the case {\it (i)}, and $\g>1,\ \al<1/2$
in the case {\it (ii)}), was proved by Heyde and Dai \cite{HD}
(see also Taniguchi and Kakizawa \cite{TK}, Theorems 6.4.1 and 6.4.2).}
\end{rem}

\subsection{Robustness to small trends of nonparametric estimation}

The next result shows that a small trend of the form $|M(t)|\leq C|t|^{-\beta}$
does not effect the asymptotic properties of the estimator of a linear spectral
functional $J(f)$, that is, the nonparametric estimation procedure is robust
to the presence of a small trend in the model.
\begin{thm}
\label{T5}
Suppose that the assumptions of Theorems \ref{TT1} and \ref{T4} are fulfilled.
Then the statistic $J(I^h_{TY})$ is consistent and asymptotically normal estimator
for functional $J(f)$ with asymptotic variance $\si^2_h(J)$ given by \eqref{tsigma}
and \eqref{eh}, that is, the asymptotic relation (\ref{t8}) is satisfied with
$I^h_{TX}(\la)$ replaced by the contaminated periodogram $I^h_{TY}(\la)$:
\bea
\label{10R}
T^{1/2}\left[J(I^h_{TY})- J(f)\right] \ConvD \eta\quad {\rm as}\quad T\to\infty,
\eea
where $\eta$ is $N(0,\si^2_h(J))$ with $\si^2_h(J)$ given by \eqref{tsigma}
and \eqref{eh}.
\end{thm}

\subsection{Robustness to small trends of parametric estimation}

The next result shows that a small trend of the form $|M(t)|\leq C|t|^{-\beta}$,
$\be>1/4$, does not effect the asymptotic properties of the Whittle estimator
of an unknown spectral parameter $\theta$, that is, the Whittle parametric
estimation procedure is robust to the presence of a small trend in the model.
\begin{thm}
\label{T5P}
Suppose that the assumptions of Theorem \ref{T4} with $g=f^{-1}(\la, \theta)\cd w(\la)$
are satisfied. Then under the conditions of Theorems \ref{TAN} the Whittle
estimator $\widehat\theta_{TY}$, constructed on the basis
of the contaminated periodogram $I_{T,Y}^h(\la)$, is consistent and asymptotically
normal estimator for an unknown spectral parameter $\theta$, that is,
the asymptotic relation (\ref{aW1}) is satisfied with $I^h_{TX}(\la)$ replaced by
the contaminated periodogram $I^h_{TY}(\la)$:
\bea
\label{W1C}
T^{1/2}\left(\widehat\theta_{TY}-\theta_0\right)\ConvD N_p\left(0, R(\theta_0)\right)
\quad {\rm as}\quad T\to\infty,
\eea
where the matrix $R(\theta_0)$ is defined in (\ref{aW2}).
\end{thm}

\begin{rem}
{\rm In the non-tapered case, Theorems \ref{T4} -- \ref{T5P} were proved in Ginovyan and Sahakyan \cite{GS2016}.
In the tapered case, the theorems can be proved similarly by
using the  tapered tools stated in Section \ref{methods}.}
\end{rem}

\section{Methods and tools}
\label{methods}

In this section we briefly discuss the methods and tools, used to prove
the results stated in Sections \ref{NP}--\ref{Rob}.

\subsection{Approximation of traces of products of Toeplitz matrices and operators.}

The trace approximation problem for truncated Toeplitz operators and matrices
has been discussed in detail in the survey paper Ginovyan et al. \cite{GST2014}
in the non-tapered case.
Here we present some important results in the tapered case, which were used
to prove the results stated in Sections \ref{NP}--\ref{Rob}.

Let $\psi(\la)$ be an integrable real symmetric function defined on
$[-\pi, \pi]$, and let $h(t)$, $t\in[0,1]$ be a taper function.
For  $T=1, 2,\ldots$, the {\it $(T\times T)$-truncated tapered Toeplitz matrix\/}
generated by $\psi$ and $h$, denoted by $B_T^h(\psi)$, is defined by the following
equation:
\beq
\label{1-7M}
B_T^h(\psi):=\|\widehat\psi(t-s)h_T(t)h_T(s)\|_{t,s=1,2\ldots,T},
\eeq
where $\widehat\psi(t)$ $(t\in \mathbb{Z})$ are the Fourier coefficients of $\psi$.

Given a real number $T>0$ and an integrable real symmetric function $\psi(\la)$
defined on $\mathbb{R}$, the {\it $T$-truncated tapered Toeplitz operator\/}
(also called {\it tapered Wiener-Hopf operator}) generated by $\psi$ and a taper function $h$,
denoted by $W_T^h(\psi)$ is defined as follows:
\begin{equation}
\label{1-7}
[{W}^h_T(\psi)u](t)=\int_0^T\hat\psi(t-s)u(s)h_T(s)ds,
\q u(s)\in L^2([0,T]; h_T),
\end{equation}
where $\hat\psi(\cdot)$ is the Fourier transform of $\psi(\cdot)$,
and $L^2([0,T]; h_T)$ denotes the weighted $L^2$-space with respect to
the measure $h_T(t)dt$.

Let $h$ be a taper function satisfying Assumption \ref{(T)},
and let $A^h_T(\psi)$ be either the $T\times T$ tapered Toeplitz matrix $B^h_T(\psi)$,
or the $T$-truncated tapered Toeplitz operator $W^h_T(\psi)$
generated by a function $\psi$ (see \eqref{1-7M} and \eqref{1-7}).

Observe that, in view of  \eqref{t4},  \eqref{t55}, \eqref{1-7M} and \eqref{1-7}, we have
%for a single matrix $B_T(f)$ we have

%C_T:= 2\pi H_{2,T}(0)=2\pi\int_0^Th_T^2(t)dt=2\pi H_2\,T
\begin{equation}
\label{MT-00}
\frac1{T} \tr\left[A^h_T(\psi)\right]=\frac1T\cd \widehat \psi(0)\cd\int_0^Th_T^2(t)dt
= 2\pi H_2\int_{\Lambda} \psi(\la)d\la.
\end{equation}
What happens to the relation (\ref{MT-00}) when $A^h_T(\psi)$
is replaced by a product of Toeplitz matrices (or operators)?
Observe that the product of Toeplitz matrices (resp. operators)
is not a Toeplitz matrix (resp. operator).

The idea is to approximate the trace of the product of
Toeplitz matrices (resp. operators) by the trace of a Toeplitz matrix
(resp. operator) generated by the product of the generating functions.
More precisely, let $\{\psi_1,\psi_2,\ldots,\psi_m\}$ be a
collection of integrable real symmetric functions defined on $\Lambda$.
Let $A_T^h(\psi_i)$ be either the $T\times T$ tapered Toeplitz matrix $B_T^h(\psi_i)$,
or the $T$-truncated tapered Toeplitz operator $W_T^h(\psi_i)$
generated by a function $\psi_i$ and a taper function $h$. Define
\begin{eqnarray}
\label{n 4-4}
\nonumber
S_{A,\mathcal{H},h}(T):=\frac1T\tr\left[\prod_{i=1}^m A_T^h(\psi_i)\right
],\q M_{\Lambda,\mathcal{H},h}:=(2\pi)^{m-1}H_m\int_{\Lambda}
\left[\prod_{i=1}^m \psi_i(\la)\right]\,d\la,
\end{eqnarray}
and let
\begin{eqnarray}
\label{n 4-5}
&&\De(T):=\De_{A,{\Lambda},\mathcal{H},h}(T)=|S_{A,\mathcal{H},h}(T)-
M_{{\Lambda},\mathcal{H},h}|.
\end{eqnarray}

\begin{pp}
\label{PT5}
Let $\De(T)$ be as in (\ref{n 4-5}). Each of the following conditions is sufficient for
\begin{equation}
\label{n 4-7}
\De(T)=o(1) \q{\rm as} \q T\to\f.
\end{equation}
\begin{itemize}
\item[{\bf(C1)}]
$\psi_i\in L^1(\Lambda)\cap L^{p_i}(\Lambda)$, $p_i>1$,
$i=1,2,\ldots,m$,
with $1/p_1+\cdots+1/p_m\le1$.

\item[{\bf(C2)}]
The function $\varphi({\bf u})$ defined by
%
%e44 ###
\begin{equation}
\label{n 4-6}
\varphi({\bf u}):\,=\,
\int_\Lambda \psi_1(\la)\psi_2(\la-u_1)\psi_3(\la-u_2)\cdots \psi_m(\la-u_{m-1})\,d\la,
\end{equation}
where ${\bf u}=(u_1,u_2,\ldots,u_{m-1})\in\Lambda^{m-1}$,
belongs to $L^{m-2}(\Lambda^{m-1})$ and is continuous at
${\bf0}=(0,0,\ldots,0)\in\Lambda^{m-1}$.
\end{itemize}
\end{pp}

\begin{rem}
\label{rem4-1}
{\rm In the non-tapered case, Proposition \ref{PT5} was proved in Ginovyan et al. \cite{GST2014},
while in the tapered case, it was proved in Ginovyan \cite{G2020g}.
Proposition \ref{PT5} was used to prove Theorems \ref{T-CLT}, \ref{TT1}, \ref{TAN} and \ref{CT1}.}
\end{rem}
\begin{rem}
{\rm More results concerning the trace approximation problem for truncated Toeplitz
operators and matrices can be found in Ginovyan and Sahakyan \cite{GS2012,GS2013}, and
in Ginovyan et al. \cite{GST2014}.}
\end{rem}

\sn{Central limit theorems for tapered quadratic functionals} % $Q_T^h$}
\label{CLT}
In this subsection we state central limit theorems for tapered quadratic
functional $Q^h_T$ given by (\ref{tq-1}), which were used to prove
the results stated in Sections \ref{NP}--\ref{Rob}.

Let $A_T^h(f)$ be either the $T\times T$ tapered Toeplitz matrix $B_T^h(f)$,
or the $T$-truncated tapered Toeplitz operator $W_T^h(f)$
generated by the spectral density $f$ and taper $h$, and let $A_T^h(g)$
denote either the $T\times T$ tapered Toeplitz matrix, or the $T$-truncated tapered
Toeplitz operator generated by the functions $g$ and $h$
(for definitions see formulas \eqref{1-7M} and \eqref{1-7}).
Similar to the non-tapered case, we have the following results
(cf. Ginovyan et al. \cite{GST2014}, Grenander and Szeg\H{o} \cite{GS}, Ibragimov \cite{I1963}).
%, see also the proof of Lemma \ref{2-9} below).

\begin{itemize}
\item[1.]
The quadratic functional $Q^h_T$ in (\ref{tq-1}) has the same distribution as the sum
$\sum_{j = 1}^\f\la_{j,T}^2\xi_j^2$, where \mbox{$\{\xi_j, j\ge1\}$} are independent
$N(0,1)$ Gaussian random variables and $\{\la_{j,T}, j\ge1\}$ are the eigenvalues of
the operator $A_T^h(f)\,A_T^h(g)$.

\item[2.]
The characteristic function $\I(t)$ of $Q^h_T$ is given by formula:
\begin{equation}
\label{MTc-05}
\I(t) = \prod_{j = 1}^\f|1 - 2it\la_{j,T}|^{-1/2}.
\end{equation}

\item[3.]
The $k$--th order cumulant $\chi_k(Q^h_T)$ of $Q^h_T$ is given by formula:
\begin{eqnarray}
\label{MTc-5}
\chi_k(Q_T) = 2^{k-1}(k-1)! \sum_{j = 1}^\f\la_{j,T}^k
=2^{k-1}(k-1)!\, \tr\,[A_T^h(f)\,A_T^h(g)]^k.
\end{eqnarray}
\end{itemize}

%Then by \eqref{MTc-5} we have
%\beq\label{1-6}
%\chi_k(\widetilde Q_T^h)= \left \{
%\begin{array}{ll}
%0, & \mbox{for $k = 1$}\\
%\\[-1mm]
% T^{-k/2}2^{k-1}(k-1)! \,
%\text{tr} \, [A_T^h(f)A_T^h(g)]^k, & \mbox{for $k \ge 2$.}
%\end{array}
%\right.
%\eeq
Thus, to describe the asymptotic distribution of the quadratic functional $Q^h_T$,
we have to control the traces and eigenvalues of the products of truncated
tapered Toeplitz operators and matrices.

\ssn{CLT for Gaussian models}
We assume that the model process $X(t)$ is Gaussian, and
with no loss of generality, that $g\ge 0$.
We will use the following notation.
By $\widetilde Q^h_T$ we denote the standard normalized
quadratic functional:
\begin{equation}
\label{7-8}
 \widetilde Q^h_T= T^{-1/2}\,\left(Q^h_T-\E [Q^h_T]\right).
\end{equation}

Also, we set
\beq\label{1-8}
\si^2_h: = 16\pi^3H_4 \int_{\Lambda} f^2(\la)g^2(\la)\,d\la,
\eeq
where $H_4$ is as in \eqref{t4}.
The notation
\begin{equation}
\label{1-55}
\widetilde Q^h_T\ConvD \eta \sim N(0,\sigma^2_h) \q {\rm as}\q T\to\f
\end{equation}
will mean that the distribution of the random variable
$\widetilde Q^h_T$ tends (as $T\to\f$) to the centered normal
distribution with variance $\sigma^2_h$.

The following theorems were  proved in
Ginovyan and Sahakyan \cite{GS2019a}.

\begin{thm}\label{th2}
Assume that $f\cdot g\in L^1(\Lambda)\cap L^2(\Lambda)$,
the taper function $h$ satisfies Assumption \ref{(T)}, % (see Section \ref{H-6.1}),
and for $T\to \infty$
\beq\label{1-11}
\chi_2(\widetilde Q^h_T)=\frac2T\tr\bigl[A^h_T(f)A^h_T(g)\bigr]^2 \longrightarrow \si^2_h,
\eeq
where $\si^2_h$ is as in \eqref{1-8}.
Then $\widetilde Q^h_T\ConvD \eta \sim N(0,\sigma^2_h)$ as $T\to\f$.
\end{thm}

\begin{thm}\label{th3}
Assume that the function
\beq\label{1-10}
\varphi(x_1, x_2,x_3)=\int_{\Lambda}
f(u)g(u-x_1)f(u-x_2)g(u-x_3)\,du
\eeq
belongs to $L^2(\Lambda^3)$ and is continuous at $(0,0,0)$,
and the taper function $h$ satisfies Assumption \ref{(T)}.
%(see Section \ref{H-6.1}).
Then $\widetilde Q^h_T\ConvD \eta \sim N(0,\sigma^2_h)$ as $T\to\f$.
\end{thm}

\begin{thm}\label{th1}
Assume that $f(\la)\in L^p(\Lambda)$ $(p\ge1)$
and  $g(\la)\in L^q(\Lambda)$ $(q\ge1)$
with $1/p+1/q\le1/2$, and the taper function $h$ satisfies Assumption \ref{(T)}.
Then $\widetilde Q^h_T\ConvD \eta \sim N(0,\sigma^2_h)$ as $T\to\f$.
\end{thm}

\begin{thm}\label{th4}
Let $f\in L^2(\Lambda)$, \,$g\in L^2(\Lambda)$, $fg\in L^2(\Lambda)$,
\begin{equation}
\label{1-12} \int_{\Lambda}
 f^2(\la)g^2(\la-\mu)\,d\la \longrightarrow
 \int_{\Lambda} f^2(\la)g^2(\la)\,d\la \quad {\rm as} \quad \mu\to0,
\end{equation}
and let the taper function $h$ satisfy Assumption \ref{(T)}.
% (see Section \ref{H-6.1}).
Then $\widetilde Q^h_T\ConvD \eta \sim N(0,\sigma^2_h)$ as $T\to\f$.
\end{thm}

To state the next theorem, we recall the class $SV_0(\mathbb{R})$ of slowly
varying functions at zero
$u(\la)$, $\la\in\mathbb{R}$, satisfying the following conditions:
for some $a>0$, $u(\la)$ is bounded on $[-a,a]$, $\lim_{\la\to0}u(\la)=0,$ \
$u(\la)=u(-\la)$ and $0<u(\la)<u(\mu)$\ for\ $0<\la<\mu<a$.

\begin{thm}\label{th5}
Assume that the functions $f$ and $g$ are integrable on $\mathbb{R}$
and bounded outside any neighborhood of the origin, and satisfy for some $a>0$
\beq \label{m-0}
f(\lambda)\le |\lambda|^{-\alpha}L_1(\lambda),
\q
|g(\lambda)|\le |\lambda|^{-\beta}L_2(\lambda),\q \lambda\in [-a,a],
\eeq
for some $\alpha<1, \ \beta<1$ with $\alpha+\beta\le1/2$, where $L_1(x)$ and
$L_2(x)$ are slowly varying functions at zero satisfying
\bea
\label{m-00}
L_i\in SV_0(\mathbb{R}), \ \ \lambda^{-(\alpha+\beta)}L_i(\lambda)\in L^2[-a,a], \ i=1,2.
\eea
Also, let the taper function $h$ satisfy Assumption \ref{(T)}.
% (see Section \ref{H-6.1}).
Then $\widetilde Q^h_T\ConvD \eta \sim N(0,\sigma^2_h)$ as $T\to\f$.
\end{thm}
The conditions $\al<1$ and $\be<1$ in Theorem \ref{th5} ensure that the
Fourier transforms of $f$ and $g$ are well defined. Observe that when $\al>0$
the process $X(t)$ may exhibit long-range dependence.
We also allow here  $\alpha+\beta$ to assume the critical value 1/2.
The assumptions $f\cdot g\in L^1(\Lambda)$, $f,g\in L^\infty(\Lambda\setminus [-a,a])$
and (\ref{m-00}) imply that $f\cdot g \in L^2(\Lambda)$, so that the variance
$\sigma^2_h$ in (\ref{1-8}) is finite.

%\paragraph{CLT for L\'evy-driven stationary linear models.}
\ssn{CLT for L\'evy-driven stationary linear models}
Now we assume that the underlying model $X(t)$ is a L\'evy-driven stationary
linear process defined by \eqref{clp}, where $a(\cdot)$ is a filter from
$L^2(\mathbb{R})$, and $\xi(t)$ is a L\'evy process satisfying the conditions:
$\E \xi(t)=0$, $\E \xi^2(1)=1$ and $\E\xi^4(1)<\infty$.

The central limit theorem that follows was proved in Ginovyan and Sahakyan \cite{GS2019}.
\begin{thm}
\label{TCLT}
Assume that the filter $a(\cdot)$
and the generating kernel $\widehat g(\cdot)$ are such that
\begin{equation}
\label{eq:CLT condition-M}
a(\cdot)\in L^p(\mathbb{R})\cap L^2(\mathbb{R}),\quad \widehat g(\cdot)\in L^q(\mathbb{R}),
\q1\le p,q\le 2,\quad \frac{2}{p}+\frac{1}{q}\ge \frac{5}{2},
\end{equation}
and the taper $h$ satisfies Assumption \ref{(T)}. % (see Section \ref{H-6.1}).
Then $\widetilde Q^h_T\ConvD \eta \sim N(0,\sigma^2_{L,h})$ as $T\to\f$, where
\beq
\label{tsigma-1}
\sigma^2_{L,h}=16\pi^3H_4\int_{\mathbb{R}} f^2(\lambda) g^2(\lambda) d\lambda
+\kappa_4 4\pi^2H_4\left[\int_{\mathbb{R}} f(\lambda) g(\lambda) d\lambda\right]^2,
\eeq
where $H_4$ is as in \eqref{t4}.
\end{thm}
\begin{rem}
{\rm
Notice that if the underlying process $X(t)$ is Gaussian,
then in formula (\ref{tsigma-1}) we have only the first term and so
$\sigma^2_{L,h}=\sigma^2_{h}$ (see \eqref{1-8}), because in this case $\kappa_4=0$.
On the other hand, the condition \eqref{eq:CLT condition-M} is more restrictive
than the conditions in Theorems \ref{th2} - \ref{th5}.
Thus, for Gaussian processes Theorems \ref{th2} - \ref{th5} improve Theorem \ref{TCLT}.
For non-tapered case Theorem \ref{TCLT}  was proved in Bai et al. \cite{BGT2}.}
\end{rem}

\subsection{Fej\'er-type kernels and singular integrals}

We define Fej\'er-type tapered kernels and singular integrals,
and state some of their properties. % that were used to prove the limit theorems stated in Section \ref{App}.

For a number $k$ ($k=2,3,\ldots$) and a taper function $h$ satisfying Assumption \ref{(T)}
consider the following Fej\'er-type tapered kernel function:
\beq \label{kerN}
F^h_{k,T}(\uu):=\frac{H_T(\uu)}{(2\pi)^{k-1}H_{k,T}(0)}, \q
\uu =(u_1, \ldots, u_{k-1})\in \mathbb{R}^{k-1},
%\begin{cases}
%\frac{H_T(\uu)}{(2\pi)^{k-1}H_{k,T}(0)},
% &  \mbox {if \, $H_{k,T}(0)\neq 0$}\\
%0, & \mbox {if \, $H_{k,T}(0)=0$},
%           \end{cases}.
\eeq
where
\beq \label{ker1N}
H_T(\uu):={H_{1,T}(u_1)\cdots H_{1,T}(u_{k-1}) H_{1,T}\left(-\sum_{j=1}^{k-1}u_j\right)},
\eeq
and the function $H_{k,T}(\cdot)$ is defined by (\ref{t3})
with $H_{k,T}(0)=T\cdot H_k\neq0$ (see (\ref{t4})).

%The proofs of lemmas that follow can be found in Ginovyan and Sahakyan \cite{GS2019}.
The next result shows that, similar to the classical Fej\'er kernel, the tapered kernel
$F^h_{k,T}(\uu)$ is an approximation identity
(see Ginovyan and Sahakyan \cite{GS2019}, Lemma 3.4).
\begin{pp}
\label{L3}
For any $k=2,3,\ldots$ and a taper function $h$ satisfying Assumption \ref{(T)}
the kernel $F^h_{k,T}(\uu)$, $\uu =(u_1, \ldots, u_{k-1})\in {\mathbb R }^{k-1}$,
possesses the following properties:
\begin{itemize}
\item [a)]
$\sup_{T>0}\int_{{\mathbb R}^{k-1}}\left |F^h_{k,T}(\uu)\right|\,d\uu =C_1<\f;$
\item [b)] $\int_{{\mathbb R}^{k-1}}F^h_{k,T}(\uu)\,d\uu =1;$
\item [c)] $\lim_{T\to\f}\int_{{\mathbb E}^c_\de}\left|F^h_{k,T}(\uu)\right|\,d\uu=0$
 for any $\de>0;$
\item [d)] If $k>2$ for any $\de>0$ there exists a constant $M_\de>0$ such that for $T>0$
\beq\label{z0}
\left\|F^h_{k,T}\right\|_{L^{p_k}({\mathbb{E}^c_\de})}\le M_\de,
\eeq
where  $p_k=\frac {k-2}{k-3}$ for $k>3$,  $p_3=\infty $ and
$$
{\mathbb E}_\de^c={\mathbb R}^{k-1}\sm{\mathbb E}_\de,\quad
{\mathbb E}_\de=\{\uu =(u_1, \ldots, u_{k-1})\in {\mathbb R }^{k-1}:
\, |u_i|\le\de, \, i=1,\ldots,k-1\}.
$$
\item [e)]
If the function $Q\in L^1(\mathbb{R}^{k-1})\bigcap L^{k-2}(\mathbb{R}^{k-1})$
and is continuous at $\vv =(v_1, \ldots, v_{k-1})$ ($L^0$ is the space of measurable functions), then
\beq\label{2-3-1}
\lim_{T\to\infty} \int_{\mathbb{R}^{k-1}}Q(\uu+\vv)F^h_{k,T}(\uu)d\uu=Q(\vv).
\eeq
\end{itemize}
\end{pp}

Denote
\beq\label{zz}
\Delta_{2,T}^h : = \int_{\mathbb{R}^2}f(\la)g(\la+\mu)F^h_{2,T}(\mu)d\la d\mu
-\int_{\mathbb{R}}f(\la)g(\la)d\la,
\eeq
where $F^h_{2,T}(\mu)$ is given by (\ref{kerN}) and (\ref{ker1N}).

The next two propositions give information on the rate of convergence
to zero of $\Delta_{2,T}^h$ as $T\to\f$ (see Ginovyan and Sahakyan \cite{GS2019}, Lemmas 4.1 and 4.2).
\begin{pp}\label{L011}
Assume that Assumptions \ref{(T)} and \ref{(A2)} are satisfied.
Then the following asymptotic relation holds: % as $T\to\f$:
\beq\label{z111}
% \int_{\mathbb{R}}\int_{\mathbb{R}}f(\la)g(\la+\mu)\Phi^h_{2,T}(\mu)d\la d\mu
%=\int_{\mathbb{R}}f(\la)g(\la)d\la + o\left(T^{-1/2}\right),
\Delta_{2,T}^h=  o\left(T^{-1/2}\right)\q {\rm as} \q T\to\f.
\eeq
\end{pp}
\begin{pp}
\label{L01}
Assume that Assumptions \ref{(T)} and \ref{(A2')} are satisfied.
Then the following inequality holds:
\beq\label{mz1}
|\Delta_{2,T}^h|\leq C_h\begin{cases}
T^{-(\be_1+\be_2)},&\text{if}\  \ \be_1+\be_2<1\\
T^{-1}\ln T,&\text{if}\  \ \be_1+\be_2=1\\
T^{-1},&\text{if} \ \ \be_1+\be_2>1,
\end{cases}\qquad T>0,
\eeq
where $C_h$ is a constant depending on $h$.
\end{pp}
Notice that for non-tapered case ($h(t)={\mathbb I}_{[0,1]}(t)$), the above stated
results were proved in Ginovyan and Sahakyan \cite{GS2007}
(see also Ginovyan and Sahakyan \cite{GS2012,GS2013}).
In the d.t.\ tapered case, Proposition \ref{L011} under different conditions
was proved in Dahlhaus \cite{D1}.

%\section{Proofs}
%\label{Proofs}

\end{document}